\documentclass[10pt]{article}

\usepackage{a4wide}
\usepackage{amssymb}
\usepackage{amsfonts}
\usepackage{amsmath}
\input xy
\xyoption{arrow}
\xyoption{matrix}

\date{}

\newtheorem{proposition}{Proposition}[section]
\newtheorem{theorem}[proposition]{Theorem}
\newtheorem{lemma}[proposition]{Lemma}

\newtheorem{corollary}[proposition]{Corollary}

\def\GK{{\rm  GK}\,}
\def\Kdim{{\rm K.dim }\,}

\def\der{\partial }

\def\nFM0{{\nu }_{F,M_0}}
\def\nFN0{{\nu }_{F,N_0}}
\def\nGN0{{\nu }_{G,N_0}}

\def\N0{ {\bf N}_0 }

\def\t{\otimes}

\def\ra{\rightarrow}

\def\lra{\leftrightarrow}
\def\Xpm{X^{\pm }}

\def\s{\sigma}
\def\Z{{\mathbb Z }}

\def\l1{{\lambda}_1}

\def\a{\alpha}
\def\a0{ {\alpha }_0}
\def\a1{ {\alpha }_1}

\def\l{\lambda}

\def\nFGM0{{\nu }_{F,G,M_0}}

\def\nFN0{{\nu}_{F,N_0}}

\def\sm{{\sigma}^m}

\def\sm1{{\sigma}^{-1}}

\def\smtp1{{\sigma}^{-t+1}}

\def\S1{S^{-1}}

\def\Xpm1{X^{\pm 1}_1}

\def\sPM1{{\sigma }^{\pm 1}}
\def\sMP1{{\sigma }^{\mp 1 }}

\def\b{\beta}
\def\d{\delta}

\def\di{{\rm d.ind}}

\def\L{\Lambda}

\def\G{\Gamma}

\def\OO{{\cal O}}
\def\CA{{\cal A}}

\def\CD{{\cal D}}

\def\Ytm1{Y^{t-1}}
\def\Yim1{Y^{i-1}}

\def\CG{{\cal G}}
\def\CH{{\cal H}}
\def\ass{{\rm ass}}

\def\i{{\bf i}}

\def\Aut{{\rm Aut}}

\def\bA{\overline{A}}
\def\Der{{\rm Der }}
\def\ad{{\rm ad }}
\def\dim{{\rm dim }}

\def\ker{ {\rm ker } }

\def\D{ \Delta }

\def\SL2Z{ {\rm SL}_2({\bf Z}) }

\def\Gp1{ G^{1 , 1 } }
\def\P11{ P^{-1 , 1 } }
\def\Pp1{ P^{1 , 1 } }

\def\nCLsr{{}^\nu\kern-2pt {\cal L}^{\sigma , \rho  }}
\def\nP{{}^\nu \kern-2pt P}
\def\nL{{}^\nu\kern-2pt L}
\def\nLL{{}^\nu\kern-2pt \Lambda}
\def\nPsr{{}^\nu\kern-2pt P^{\sigma , \rho  }}
\def\nLsr{{}^\nu\kern-2pt L^{\sigma , \rho  }}
\def\nuCL{{}^\nu\kern-2pt  {\cal L}}
\def\nCLsr{{}^\nu\kern-2pt {\cal L}^{\sigma , \rho  }}
\def\nCL1m{{}^\nu\kern-2pt {\cal L}^{-1 , 1  }}
\def\x1nu{x^\frac{1}{\nu}}
\def\xm1nu{x^{-\frac{1}{\nu}}}

\def\ra{\rightarrow }

\def\CB{{\cal B}}

\def\CI{{\cal I}}

\def\CH{ {\cal H}}
\def\CP{ {\cal P}}

\def\nAM0{{\nu }_{{\cal A},M_0}}
\def\nAN0{{\nu }_{{\cal A},N_0}}

\def\Kdim{ {\rm Kdim } }

\def\End{ {\rm End }}
\def\Der{ {\rm Der }}

\def\CP{ {\cal P }}
\def\det{ {\rm det }}
\def\ad{ {\rm ad }}

\def\bp{\overline{p}}

\def\ga{\mathfrak{a}}
\def\gb{\mathfrak{b}}
\def\gc{\mathfrak{c}}

\def\gn{\mathfrak{n}}
\def\gm{\mathfrak{m}}
\def\gp{\mathfrak{p}}

\def\bJ{\overline{J}}

\def\j{{\bf j}}

\def\Spec{{\rm Spec}}

\def\II{{\bf I}}
\def\JJ{{\bf J}}

\def\bJ{\overline{J}}

\def\SL{{\rm SL}}
\def\j{{\bf j}}

\def\Spec{{\rm Spec}}

\def\II{{\bf I}}
\def\JJ{{\bf J}}

\def\di!{\frac{\der^i}{i!}}
\def\dik!{\frac{\der^k_i}{k!}}

\def\gl{\mathfrak{l}}

\def\N{\mathbb{N}}

\def\0{\overline{0}}
\def\1{\overline{1}}

\def\Ln1{\L_{n,\overline{1}}}

\def\a1{a_{\overline{1}}}

\def\S{\Sigma}

\def\CU{{\cal U}}

\def\vn1{\overrightarrow{n-1}}

\def\gl{{\rm gl}}
\def\sl{{\rm sl}}

\def\Sub{{\rm Sub}}

\def\mJ{\mathbb{J}}
\def\mI{\mathbb{I}}

\def\K1{{\rm K}_1}

\def\hmI1{\widehat{\mI_1}}
\def\tmI1{\widetilde{\mI_1}}
\def\tmJ1{\widetilde{\mJ_1}}
\def\hB1{\widehat{B_1}}
\def\hCB1{\widehat{\CB_1}}

\def\bS{\overline{S}}

\def\Den{{\rm Den}}

\def\Den{{\rm Den}}

\def\ga{\mathfrak{a}}

\def\tor{{\rm tor}}

\def\bE{\overline{E}}

\def\sl2{\mathfrak{sl}_2}
\def\gl2{\mathfrak{gl}_2}
\def\gs{\mathbb{s}}
\def\gs{\mathfrak{s}}

\def\Irr{{\rm Irr}}
\def\b1{\overline{1}}

\begin{document}

\author{V. V. \  Bavula   
}

\title{$\D$-locally nilpotent algebras, their ideal structure and simplicity criteria}

\maketitle
\begin{abstract}

The class of $\D$-locally nilpotent algebras (introduced in the paper) is a wide generalization of the algebras of differential operators on commutative algebras. Examples includes all the rings $\CD (A)$ of differential operators on commutative algebras (in arbitrary characteristic), all subalgebras of $\CD (A)$ that contain the algebra $A$, the universal enveloping algebras of nilpotent, solvable and semi-simple Lie algebras, the Poisson universal enveloping algebra of an arbitrary Poisson algebra, iterated Ore extensions $A[x_1, \ldots , x_n ; \d_1 , \ldots , \d_n]$, certain generalized Weyl algebras,  and others.

 In \cite{SimCrit-difop}, simplicity criteria are given  for the algebras differential operators on commutative algebras (it was a long standing problem). The aim of the paper is to describe the ideal structure of $\D$-locally nilpotent algebras and as a corollary to give simplicity criteria for them (it is a generalization of the results of \cite{SimCrit-difop}). Examples are considered. \\
 


$${\bf Contents}$$
\begin{enumerate}
\item Introduction. 
\item Simplicity criteria of $\D$-locally nilpotent algebras: Proofs.
\item Simplicity criteria for subalgebras of $\CD (\CA )$ that contain $\CA$.
\end{enumerate}

\end{abstract}


\section{Introduction}

In the paper,  $K$ is a field of arbitrary  characteristic (not necessarily algebraically closed); module means a left
module; for a commutative algebra $A$, $\CD (A)$ is the algebra of differential operators on $A$ and $\Der_K(A)$ is the left $A$-module of $K$-derivations of $A$.\\

{\bf Simplicity criteria for the algebra $\CD (\CA )$ of differential operators  on  the algebra $\CA$ which is a domain of essentially finite type.} Theorem \ref{CX9Jul19} and Theorem \ref{CXA9Jul19} are simplicity criteria for the algebra $\CD (\CA )$ where $\CA$ is a domain of essentially finite type over a perfect field (Theorem \ref{CX9Jul19}) and  a commutative algebra over an arbitrary field (Theorem \ref{CXA9Jul19}),  respectively. 

The aim of the paper is to generalize the above results for a large class of algebras -- the $\D$-{\em locally nilpotent algebras} -- which includes  the algebra $\CD (A)$ of differential operators on a commutative algebra $A$ and all its  subalgebras that contain the algebra $A$. The last class of algebras contains many exotic algebras (non-Noetherial and not finitely generated).

\begin{theorem}\label{CX9Jul19}
(\cite[Theorem 1.1]{SimCrit-difop}) Let  a $K$-algebra $\CA$ be a commutative domain of essentially finite type over a perfect field $K$ and $\ga_r$ be its Jacobial ideal. The following statements are equivalent:
\begin{enumerate}
\item The algebra $\CD (\CA )$ of differential operators on $\CA$ is a simple algebra.
\item For all $i\geq 1$, $\CD (\CA ) \ga_r^i\CD (\CA )= \CD (\CA )$. 
\item For all $k\geq 1$, $\i\in \II_r$ and $\j\in \JJ_r$,  $\CD (\CA ) \D (\i, \j )^k\CD (\CA )= \CD (\CA )$.
\end{enumerate}
\end{theorem}
The elements $\D (\i, \j )$ are defined in Section \ref{SIMCRIT} (they are  the non-zero minors of maximal rank in the Jacobi matrix). Theorem \ref{CX9Jul19} presents a  short proof of an important result in the area of differential operators that if the algebra $\CA$ is a smooth then the algebra $\CD (\CA)$ is simple: If the algebra $\CA$ is smooth, i.e. $\ga_r=\CA$ (the Jacobian Criterion of Regularity), then  by the second condition of Theorem \ref{CX9Jul19} the algebra $\CD (\CA)$ is a simple algebra. 
 Theorem \ref{CX9Jul19} reveals the reason why for  some singular algebras $\CA$ their rings of differential operators are simple algebras.  For example, this is the case for the cusp. 


\begin{theorem}\label{1Jan20}
 Let  a $K$-algebra $\CA$ be a commutative domain of essentially finite type over a perfect field $K$  and $\ga_r$ be its Jacobian ideal. The following statements are equivalent:
\begin{enumerate}
\item The algebra $\CD (\CA )$ of differential operators on $\CA$ is a simple algebra.
\item For every maximal ideal $\gm$ of $\CA$ that contains the Jacobian ideal $\ga_r$,
the algebra $\CD (\CA)_\gm$ is a simple algebra.  
\end{enumerate}
\end{theorem}
The proof of Theorem \ref{1Jan20} is given in Section \ref{SIMCRIT}.\\

{\bf Simplicity criterion for the algebra $\CD (R )$  of differential operators  on  an arbitrary commutative algebra $R$.} An ideal $\ga $ of the algebra $R$ is called $\Der_K(R)$-{\em stable} if $\d (\ga ) \subseteq \ga$ for all $\d \in \Der_K(R)$. Theorem \ref{CXA9Jul19}.(2) is a simplicity criterion for the algebra $\CD (R )$ where $R$ is an arbitrary commutative algebra. Theorem \ref{CXA9Jul19}.(1) shows that every nonzero ideal of the algebra $\CD (R)$ meets the subalgebra $R$ of $\CD (R)$. If, in addition, the algebra $R=\CA $ is a domain of essentially finite type,  Theorem \ref{CXA9Jul19}.(3) shows that every nonzero ideal of the algebra $\CD (R)$ contains a power of the Jacobian ideal of $\CA$.

\begin{theorem}\label{CXA9Jul19}
(\cite[Theorem 1.2]{SimCrit-difop}) Let $R$ be a commutative algebra over an arbitrary field $K$. 
\begin{enumerate}
\item Let $I$ be a nonzero ideal the algebra $\CD (R)$. Then the ideal $I_0:= I\cap R$ is a nonzero $\Der_K(R)$-stable ideal of the algebra $R$ such that $\CD (R)I_0\CD (R)\cap R = I_0$. In particular, every nonzero ideal of the algebra $\CD (R)$ has nonzero intersection with $R$.
\item The ring $\CD (R)$ is not simple iff there is a proper $\Der_K(R)$-stable ideal $\ga $ of $R$ such that $\CD (R)\ga \CD (R)\cap R = \ga$.
\item Suppose, in addition, that $K$ is a perfect field and  the algebra $\CA = R$ is a domain of essentially finite type, $\ga_r$ be its Jacobian ideal, $I$ be a nonzero ideal of $\CD (\CA )$, and $I_0=I\cap \CA$. Then $\ga_r^i\subseteq I_0$ for some $i\geq 1$. 
\end{enumerate}
\end{theorem}

{\bf The $\D$-locally nilpotent modules.} The following notations will remain fixed in the paper:  $A$ is a $K$-algebra, $M$ is an $A$-module, $\emptyset \neq \D\subseteq \End_A(M)$ and $\D^i=\D\cdots \D =\{ \d_1\cdots \d_i\, | \,\d_1, \ldots , \d_i\in \D\}$ ($i\geq 1$ times), $$N(M)=N_\D (M):=\bigcup_{i\geq 0} N_\D (M)_i\;\; {\rm where}\;\; N(M)_i=N_\D (M)_i:={\rm ann}_M(\D^{i+1})=\{ m\in M\, | \, \D^{i+1}m=0\}$$ and $N(M)_{-1}:=0$. Clearly, $N(M)_{-1}\subseteq N(M)_0\subseteq \cdots \subseteq N(M)_n\subseteq \cdots $ is an ascending chain of $A$-submodules of $M$ such that 
$$\D N(M)_i\subseteq N(M)_{i-1} \;\; {\rm for\; all}\;\; i\geq 0.$$
{\em Definition.} The $A$-module $N_\D (M)$ is called the $\D$-{\bf locally nilpotent $A$-submodule} of $M$. The $A$-module $M$ is called the $\D$-{\bf locally nilpotent $A$-module} if $M=N_\D (M)$.\\

In general situation, the $A$-submodule $N_\D (M)$ of $M$ is a $\D'$-locally nilpotent $A$-module where $\D'=\{ \d'\, | \, \d \in \D\}$ and $\d'$ is the restriction of the $A$-homomorphism $\d$ to $N_\D (M)$. Abusing the language, we call the $A$-module $N_\D (M)$ the $\D$-locally nilpotent $A$-module. 

A map $f\in \End_A(M)$ is called a {\bf locally nilpotent map} if $M=\bigcup_{i\geq 0}\ker_M(f^{i+1})$. If $M$ is a $\D$-locally nilpotent $A$-module then every map $\d \in \D$ is a locally nilpotent map but {\em not}  vice versa, in general, see the example below. \\

{\em Example.} Let $M=\bigoplus_{i=1}^n Ae_i$ be a free $A$-module of rank $n\geq 2$ where the set $\{ e_1, \ldots , e_n\}$ is a free basis for $M$; $\D = \{ \d_+, \d_-\}\subseteq \End_A(M)$ where $\d_\pm (e_i)=e_{i\pm 1}$ for $i=1,\ldots , n$ and $e_{0}=e_{n+1}=0$. Clearly, $\d_\pm^n=0$, the maps $E_+=\d_+\d_-$ and $E_-=\d_-\d_+$ are nonzero  idempotents such that $E_+(e_i)=e_i$ for $i=2, \ldots , n$ and $E_+(e_1)=0$, 
 $E_-(e_i)=e_i$ for $i=1, \ldots , n-1$ and $E_-(e_n)=0$. Therefore, the $A$-module $M$ is not $\D$-locally nilpotent. In fact, $N_\D (M)=0$. 
 
 \begin{lemma}\label{b1Jan20}
If $\D$ is a finite set of commuting $A$-homomorphism of an $A$-module $M$. Then the $A$-module is $\D$-locally nilpotent iff all the maps in $\D$ are locally nilpotent maps.
\end{lemma}

{\it Proof}. Straightforward. $\Box $\\

{\bf The $\D$-locally nilpotent algebras $N_\D (E)$ where $\D\subseteq \Der_A(E)$.} Suppose, in addition, that $A$ is a subalgebra of an {\em algebra} $E$ and $\D\subseteq \Der_A(E)$, the set of $A$-{\em derivations} of the algebra $E$ ($\d\in \Der_A(E)$ if $\d$ is a derivation of the algebra $E$ and an $A$-homomorphism; in particular, $\d (A)=0$). Then $E^\D :=\bigcap_{\d \in \D} \ker_E(\d )$ is the {\bf algebra of $\D$-constants}, and $A\subseteq E^\D$.

\begin{proposition}\label{A23Dec19}
Let $A$ be a subalgebra of an algebra $E$ and $\D \subseteq \Der_A(E)$. Then
\begin{enumerate}
\item The $A$-module $N_\D (E)=\bigcup_{i\geq 0}N_\D (E)_i$ is a subalgebra of $E$ such that $A\subseteq N_\D (E)_0=E^\D$, $N_\D (E)_iN_\D (E)_j\subseteq N_\D (E)_{i+j}$ for all $i,j\geq 0$, i.e. the set $\{ N_\D (E)_i\}_{i\geq 0}$ is   an ascending filtration of the algebra $N_\D (E)$ elements of which are $A$-modules.
\item For all $i\geq 0$, $\D N_\D (E)_i\subseteq N_\D (E)_{i-1}$.
\end{enumerate}
\end{proposition}

{\em Definition.} The algebra $N_\D (E)$ is called the $\D$-{\bf locally nilpotent algebra} and the filtration $\{ N_\D (E)_i\}_{i\geq 0}$ is called the {\bf order filtration} on the algebra $N_\D (E)$. We say that an element $a\in N_\D (E)_i\backslash 
N_\D (E)_{i-1}$ has {\bf order} $i$ which is denoted by ${\rm ord} (a)=i$.\\

The $\D$-locally nilpotent algebras are the main object of study of the paper. We clarify their ideal structure and give several simplicity criteria for them. Below are examples of several large classes of $\D$-locally nilpotent algebras.\\

{\em Example {\sc (The algebras of differential operators)}.} Let $A$ be a {\em commutative} $K$-algebra, $E:=\End_K(A)\supseteq \End_A(A)\simeq A$ and $\D = \{ \ad_a\, | \, a\in A\}$ where $\ad_a: E\ra E$, $ f\mapsto [a,f]:=af-fa$ is the {\bf inner derivation} of the algebra $E$ determined by the element $a$. By the very definition, $\D \subseteq \Der_A(E)$ and 
\begin{equation}\label{NDE=CDA}
N_\D (E)=\CD (A)
\end{equation}
is the {\bf algebra of differential operators} on the algebra $A$ and the filtration $\{ N_\D (E)=\CD (A)_i\}_{i\geq 0}$ is the {\bf order filtration} on the algebra $\CD (A)$, see Section \ref{SIMCRIT} for details.\\ 
 
{\em Example {\sc (Subalgebras  of differential operators $\CD (A)$ that contain $A$)}.}   A subalgebra $R$ of the algebra $\CD (A)$ of differential operators on a commutative algebra $A$ that  contains the algebra $A$ is a $\D$-locally nilpotent algebra w.r.t. $\D =\{ \ad_a\, | \, a\in A\}$ and the induced filtration $\{R_i:=R\cap \CD (A)_i\}_{i\geq 0}$ is the $\D$-order filtration on $R$. See Proposition \ref{A25Dec19} for examples.\\
 

 {\em Example.} Suppose that the algebra $E$ admits  a set of generators $\{ a_i\, | \, i\in I\}$ such that the algebra $E$ is a $\D$-locally nilpotent algebra where $\D =\{ \ad_{a_i}\, | \, i\in I\}\subseteq \Der_{A}(E)$ and $A=Z(E)$ is the centre of the algebra $E$.  \\ 

{\em Example.} The {\bf Weyl algebra} 
$$A_n=K\langle x_1, \ldots , x_n \der_1, \ldots , \der_n\, | \, [\der_i,x_j]=\d_{ij}, \;\; x_ix_j=x_jx_i, \;\;\der_i\der_j=\der_j\der_i, \;\;1\leq i,j\leq n \rangle$$ is a $\D$-locally nilpotent algebra where $\D =\{ \ad_{x_i}, \ad_{\der_i}\, | \, i=1, \ldots , n\}$ where $[a,b]:=ab-ba$ and  $\d_{ij}$  is the Kronecker delta. If the field $K$ has characteristic zero then $Z(A_n)=K$ and the $\D$-order filtration $\{ A_{n,i}\}_{i\geq 0}$ coincides with the standard filtration on $A_n$ with respect to the canonical generators $ x_1, \ldots , x_n,  \der_1, \ldots , \der_n$ of the algebra $A_n$. Namely, $$A_{n,i}=\sum_{|\alpha |+|\beta |\leq i }Kx^\alpha \der^\beta$$ where $\alpha , \beta \in \N^n$ and $|\alpha |=\alpha_1+\cdots +\alpha_n$ for $\alpha =(\alpha_1,\ldots ,\alpha_n)$. For all $i\geq 0$, $\dim_K(A_{n,i})<\infty$.  The Weyl algebra $A_n$ is a simple algebra if the field  $K$ has characteristic zero. The  Weyl algebra  $A_n$ is also  a $\D$-locally nilpotent algebra where $\D =\{ \ad_{x_i}\, | \, i=1, \ldots , n\}$ or  $\D =\{  \ad_{\der_i}\, | \, i=1, \ldots , n\}$ but in these cases the components of the $\D$-order filtrations  are infinite dimensional. \\

{\em Example.} Let $\gn$ be a nilpotent Lie algebra. Then its  universal enveloping algebra $U=U(\gn )$ is a $\D$-locally nilpotent algebra where $\D =\{ \ad_x \, | \, x\in \gn\}\subseteq \Der_A(U)$ and $A=Z(U)=U^\D$ is the centre of the algebra $U$. \\

 {\em Example.} Let $\gs$  be a solvable Lie algebra. Then $\gn = [\gs , \gs ]$ is  a nilpotent Lie algebra and  the  universal enveloping algebra $U=U(\gs )$ is a $\D$-locally nilpotent algebra where $\D =\{ \ad_x \, | \, x\in \gn\}\subseteq \Der_A (U)$ and $A=U^\D$ is the centralizer of $\gn$ in the algebra $U$. \\
 
 {\em Example.} Let $\CG$ be a semi-simple  Lie algebra and $\CG = \gn_-\oplus \CH \oplus \gn_+$ be its triangular decomposition where $\CH$ is  a Cartan subalgebra of $\CG$.  Then the  universal enveloping algebra $U=U(\CG )$ is a $\D_+$-locally nilpotent  (resp., $\D_-$-locally nilpotent) algebra where $\D_+ =\{ \ad_x \, | \, x\in \gn_+\}\subseteq \Der_{A_+}(U)$ (resp., $\D_- =\{ \ad_x \, | \, x\in \gn_-\}\subseteq \Der_{A_-}(U)$) and $A_+=U^{\D_+}$ (resp., $A_-=U^{\D_-}$). \\
 
 {\em Example {\sc (The Poisson universal enveloping algebra of a Poisson algebra)}.} Let $\CP$ be a Poisson  algebra. In \cite{GenDefRel-PUEA}, for each Poisson algebra $\CP$ explicit sets of (associative) algebra generators and defining relations are given for  its Poisson universal enveloping algebra $\CU (\CP )$. It follows at once from this result that the  {\bf Poisson universal enveloping algebra} $\CU (\CP )$ is  a $\D$-locally nilpotent algebra w.r.t. $\D =\{ \ad_a\, | \, a\in \CP \}$. \\

 {\em Example {\sc (The algebra of Poisson differential operators)}.} Let $(\CP, \{ \cdot, \cdot \} )$ be a Poisson  algebra. In \cite{GenDefRel-PUEA}, the {\bf algebra of of Poisson differential operators} $\CP \CD (\CP)$ was introduced and studied. Since  $\CP \subseteq \CP \CD (\CP)\subseteq \CD (\CP)$, the algebra  $\CP \CD (\CP)$ is a $\D$-nilpotent algebra where $\D = \{ \ad _p\, | \, p\in \CP\}$. \\

 {\em Example.} Let $A$ be a {\em commutative} algebra and $\{\d_1, \ldots , \d_n\}$ be a set of commuting $K$-derivations of the algebra $A$, $R=A[x_1, \ldots , x_n; \d_1, \ldots , \d_n]$ be an {\em iterated Ore extension}. The algebra $R$ is generated by the algebra $A$ and elements $x_1, \ldots , x_n$ subject to the defining relations:
 $$x_ix_j=x_jx_i  \;\;(i\neq j)\;\; {\rm and} \;\; x_ia=ax_i+\d_i(a)\;\; (a\in A, \; 1\leq i \leq n).$$
The algebra $R=\bigoplus_{\alpha \in \N^n} A\d^\alpha =  \bigoplus_{\alpha \in \N^n} \d^\alpha A$ is a free left and right $A$-module with free basis $\{ \d^\alpha \}_{\alpha \in \N^n}$ where $\d^\alpha = \d_1^{\alpha_1}\cdots \d_n^{\alpha_n}$ and $\alpha = (\alpha_1 , \ldots , \alpha_n)$. In particular, the algebra $A$ is a subalgebra of $R$. The algebra $R$ is a $\D$-locally nilpotent algebra where $\D =\{ \ad_a\, | \, a\in A\}\subseteq \Der_A(R)$. Let $\{ R_i\}$ be the $\D$-order filtration on $R$. Then $R_0=R^\D=C_R(A)$ is the centralizer of the algebra $A$ in $R$. In particular, $A\subseteq R_0$ and $\sum_{i=1}^n Ax_i= \sum_{i=1}^n x_iA\subseteq \sum_{i=1}^n R_0x_i=\sum_{i=1}^n x_iR_0$.\\

 {\em Example.} Certain generalized Weyl algebras of rank $n$ are $\D$-locally nilpotent algebras, see Lemma \ref{a1Jan20} for details.\\

 {\bf The ideal structure of the $\D$-stable ideals of the $\D$-locally nilpotent algebras.} An ideal $\ga$ of an algebra $E$ is called $\D$-{\bf stable} if $\D \ga \subseteq \ga$ ($\d (\ga ) \subseteq \ga$ for all $\d \in \D$) where $\D\subseteq \Der_A(E)$. Theorem \ref{B23Dec19} describes the ideal structure of  the $\D$-stable ideals of the $\D$-locally nilpotent algebras. This theorem is used in many proofs of the paper. 
 
 \begin{theorem}\label{B23Dec19}
Let $A$ be a subalgebra of an algebra $E$, $\D\subseteq \Der_A(E)$, $\ga$ be a nonzero $\D$-stable ideal of the algebra $N_\D (E)$,  and $\ga_0:=\ga\cap N_\D (E)_0=\ga\cap E^\D$. Then
\begin{enumerate}
\item The ideal $\ga_0\neq 0$ is a nonzero ideal of the algebra $N_\D (E)_0$,  the ideal  $\ga'= N_\D (E)\ga_0 N_\D (E)$ is a nonzero ideal of the algebra $N_\D (E)$  such that $\ga'\subseteq \ga$ and $\ga'\cap N_\D (E)_0=\ga_0$. 
\item If, in addition, the algebra $N_\D (E)_0$ is a commutative algebra then $[N_\D (E)_1, \ga_0]\subseteq \ga_0$. 
\end{enumerate}
\end{theorem}

{\bf The ideal structure of  the $\D$-locally nilpotent algebras $N_\D (E)$ where $A$ is a commutative algebra.} For a subset $A'$ of the algebra $A$, $C_A(A'):=\{ a\in A\, | \, aa'=a'a$ for all $a\in A\}$ is the {\em centralizer} of the set $A'$ in $A$. The centralizer $C_A(A')$  is a subalgebra of $A$. 

 \begin{theorem}\label{23Dec19}
Let $A$ be a commutative subalgebra of an algebra $E$, $\D =\{ \ad_a\, | \, a\in A'\}$ where $A'$ is a non-empty subset of $A$ (eg, $A'=A$), and $\ga$ be a nonzero ideal of the algebra  $N_\D (E)$ and  $\ga_0:=\ga\cap N_\D (E)_0=\ga\cap E^\D$. Then
\begin{enumerate}
\item  $\ga_0\neq 0$ is a nonzero ideal of the algebra $N_\D (E)_0=E^\D=C_E(A')$. The ideal $\ga'=N_\D (E)\ga_0 N_\D (E)$ of $N_\D (E)$ is a nonzero ideal such that $\ga'\cap N_\D (E)=\ga_0$.
\item If, in addition, the algebra $N_\D (E)_0$ is a commutative algebra then $[N_\D (E)_1, \ga_0]\subseteq \ga_0$. 
\end{enumerate}
\end{theorem}

Using Theorem \ref{23Dec19}, we clarify the ideal structure of  subalgebras of $\CD (A)$ that contain $A$, Thereom \ref{C23Dec19}. This result is a generalization of Theorem \ref{CXA9Jul19}.

\begin{theorem}\label{C23Dec19}
Let $A$ be a commutative algebra, $\CD (A)$ be the algebra of differential operators on $A$, $R$ be a subalgebra  of $\CD (A)$ such that $A\subseteq R$ (eg, $R= \CD (A)$), and $R_i=R\cap \CD (A)_i$ for $i\geq 0$. If  $\ga$ is a nonzero ideal of the algebra $R$ then $\ga_0:=\ga\cap A\neq 0$ is a nonzero ideal of the algebra $A$ such that $[R_1, \ga_0]\subseteq \ga_0$ and $R\ga_0R\cap A=\ga_0$. The condition that  $[R_1, \ga_0]\subseteq \ga_0$ is equivalent to the condition that  $[D_R, \ga_0]\subseteq \ga_0$ where $D_R:=R_1\cap \Der_K(A)$ (if $R=\CD (A)$ then $D_R=\Der_K(A)$).
\end{theorem}

The proofs of Theorem \ref{B23Dec19}, Theorem \ref{23Dec19}  and Theorem \ref{C23Dec19} are given in Section \ref{PROOFS}.

\begin{corollary}\label{a5Jan20}
Let $A$ be an algebra that admits  a set of generators $\{ a_i\, | \, i\in I\}$ such that the algebra $A$ is a $\D$-locally nilpotent algebra where $\D =\{ \ad_{a_i}\, | \, i\in I\}\subseteq \Der_{Z(A)}(A)$ (eg, the universal enveloping algebra $U(\gn )$ of a nilpotent algebra). Then every nonzero ideal of $A$ meets the centre $Z(A)$ of $A$. 
\end{corollary}

{\it Proof}. The zero term of the $\D$-order filtration of the algebra $A$ is $A^\D = Z(A)$ the centre of the algebra $A$. Now, the statement follows from Theorem \ref{23Dec19}.(1). $\Box $\\

{\bf Simplicity criteria for subalgebras $R$ of $\CD (\CA )$ that contains $\CA$.} The algebras $R$ are a very wide class of algebras and they are important examples of $\D$-locally nilpotent algebras. Even in the case of the polynomial algebra $\CA = K[x]$ over a field of characteristic zero, the structure of the algebras $R\subseteq \CD (K[x])$  are not yet completely understood as there are exotic algebras (not finitely generated and not Noetherian). Some examples of such algebras are considered in Proposition \ref{A25Dec19}. In particular,  their prime spectra are classified. A submodule of a module is called an {\em essential} submodule if it meets all the nonzero submodules. 

\begin{theorem}\label{27Dec19}
Let $\CA$ be a commutative domain of essentially finite type over a perfect field, $\ga_r$ be its Jacobian ideal,  $\CD (\CA )$ be the algebra of differential operators on $\CA$, $R$ be a subalgebra of $\CD (\CA )$ that contains $\CA$ and is an essential $\CA$-submodule of $\CD (\CA )$, and $R_i=R\cap \CD (\CA )_i$ for $i\geq 0$. For each $i\geq 1$, let $\gb_i:={\rm l.ann}_\CA (\CD (\CA )_i/R_i)$ and $\gc_i:={\rm r.ann}_\CA (\CD (\CA )_i/R_i)$. 
 Then the following statements are equivalent: 
\begin{enumerate}
\item The algebra $R$ is a simple algebra.
\item For all integers $i\geq 1$, $\CD (\CA ) \ga_r^i\CD (\CA ) = \CD (\CA )$ and $R\gb_i\gc_iR=R$. 
\item The algebra $\CD (\CA )$ a simple algebra and $R\gb_i\gc_iR=R$ for all $i\geq 1$. 
\item For all integers $i\geq 1$, $\CD (\CA ) \ga_r^i\CD (\CA ) = \CD (\CA )$, $R\gb_1^2R=R$  and $R\gb_1\cdots \gb_{i-1}\gb_{i}^2R=R$. 
\item The algebra $\CD (\CA )$ a simple algebra, $R\gb_1^2R=R$  and $R\gb_1\cdots \gb_{i-1}\gb_{i}^2R=R$ for all $i\geq 2$.
\item For all integers $i\geq 1$, $\CD (\CA ) \ga_r^i\CD (\CA ) = \CD (\CA )$, $R\gc_1^2R=R$  and $R\gc_1\cdots \gc_{i-1}\gc_{i}^2R=R$. 
\item The algebra $\CD (\CA )$ a simple algebra, $R\gc_1^2R=R$  and $R\gc_1\cdots \gc_{i-1}\gc_{i}^2R=R$ for all $i\geq 2$.
\end{enumerate}
\end{theorem}

The proof of Theorem \ref{27Dec19} is given in Section \ref{SIMCRIT}. Using Theorem \ref{27Dec19}, we obtain Theorem \ref{28Dec19} which is another simplicity criterion for the algebra $R$.

\begin{theorem}\label{28Dec19}
Let $\CA$ be a commutative domain of essentially finite type over a perfect field, $\ga_r$ be its Jacobian ideal,  $\CD (\CA )$ be the algebra of differential operators on $\CA$, $R$ be a subalgebra of $\CD (\CA )$ that contains $\CA$, and $S^{-1}R=S^{-1}\CD (\CA )$ for some multiplicative subset $S$ of $\CA$. Fix elements $s_i, t_i\in S$ such that $s_i\in \gb_i$ and $t_i\in \gc_i$ for $i\geq 1$ (see Theorem \ref{27Dec19} for the definition of the ideals $\gb_i$ and $\gc_i$)
 Then the following statements are equivalent:
\begin{enumerate}
\item The algebra $R$ is a simple algebra.
\item For all integers $i\geq 1$, $\CD (\CA ) \ga_r^i\CD (\CA ) = \CD (\CA )$ and $Rs_it_iR=R$. 
\item The algebra $\CD (\CA )$ a simple algebra and $Rs_it_iR=R$ for all $i\geq 1$. 
\item For all integers $i\geq 1$, $\CD (\CA ) \ga_r^i\CD (\CA ) = \CD (\CA )$, $Rs_1^2R=R$  and $Rs_1\cdots s_{i-1}s_{i}^2R=R$. 
\item The algebra $\CD (\CA )$ a simple algebra, $Rs_1^2R=R$  and $Rs_1\cdots s_{i-1}s_{i}^2R=R$ for all $i\geq 2$.
\item For all integers $i\geq 1$, $\CD (\CA ) \ga_r^i\CD (\CA ) = \CD (\CA )$, $Rt_1^2R=R$  and $Rt_1\cdots t_{i-1}t_{i}^2R=R$. 
\item The algebra $\CD (\CA )$ a simple algebra, $Rt_1^2R=R$  and $Rt_1\cdots t_{i-1}t_{i}^2R=R$ for all $i\geq 2$.
\end{enumerate}
\end{theorem}

The proof of Theorem \ref{28Dec19} is given in Section \ref{SIMCRIT}.

\begin{theorem}\label{8Jan20}
Let $\CA$ be a commutative domain of essentially finite type over a perfect field, $\ga_r$ be its Jacobian ideal,  $\CD (\CA )$ be the algebra of differential operators on $\CA$, $R$ be a subalgebra of $\CD (\CA )$ that contains $\CA$ and is an essential $\CA$-submodule of $\CD (\CA )$, and $R_i=R\cap \CD (\CA )_i$ for $i\geq 0$. For each $i\geq 1$, let $\gb_i={\rm l.ann}_\CA (\CD (\CA )_i/R_i)$ and $\gc_i={\rm r.ann}_\CA (\CD (\CA )_i/R_i)$. 
 Then the following statements are equivalent: 
\begin{enumerate}
\item The algebra $R$ is a simple algebra.
\item The algebra $\CD (\CA )$ is a simple algebra and for every maximal ideal $\gm$ of $\CA$ that contains one of the  ideals $\gb_i$ ($i\geq 1$),
the algebra $R_\gm$ is a simple algebra.
\item For every maximal ideal $\gn$ of $\CA$ that contains the Jacobian ideal $\ga_r$,
the algebra $\CD (\CA )_{\gn}$ is a simple algebra and for every maximal ideal $\gm$ of $\CA$ that contains one of the  ideals $\gb_i$ ($i\geq 1$), the algebra $R_\gm$ is a simple algebra.
\item The algebra $\CD (\CA )$ is a simple algebra and for every maximal ideal $\gm$ of $\CA$ that contains one of the  ideals $\gc_i$ ($i\geq 1$),
the algebra $R_\gm$ is a simple algebra.
\item For every maximal ideal $\gn$ of $\CA$ that contains the Jacobian ideal $\ga_r$,
the algebra $\CD (\CA )_{\gn}$ is a simple algebra and for every maximal ideal $\gm$ of $\CA$ that contains one of the  ideals $\gc_i$ ($i\geq 1$), the algebra $R_\gm$ is a simple algebra.
\end{enumerate}
\end{theorem}

The proof of Theorem \ref{8Jan20} is given in Section \ref{SIMCRIT}.


\section{Simplicity criteria of $\D$-locally nilpotent algebras: Proofs}\label{PROOFS}

In this section proofs of Theorem \ref{A23Dec19}, Theorem \ref{B23Dec19}, Theorem \ref{23Dec19} and Theorem \ref{C23Dec19} are given. A class of generalized Weyl algebras that are $\D$-locally nilpotent algebras are considered (Lemma \ref{a1Jan20}). In  Proposition \ref{A25Dec19}, the prime spectra of  three $\D$-locally nilpotent subalgebras of the Weyl algebra $A_1$ are described. Proposition \ref{A2Jan20} shows that  $N_\D (S^{-1}E)\simeq S^{-1}N_\D (E)$ for all  regular left Ore sets $S$ that are contained in the zero component $N_\D (E)_0$ of the $\D$-order filtration. Similar results hold in more general situation but under additional condition (Proposition \ref{B2Jan20}). We study properties of denominator sets  that are generated by $\ad$-locally nilpotent elements (Proposition \ref{C2Jan20}). In particular, localizations at such denominator sets respect ideals (which is not true for arbitrary localization).

If $\d$ is a derivation of an algebra $E$ then for all elements $a,b\in E$,
\begin{equation}\label{dnab}
\d^n(ab)=\sum_{i=0}^n {n\choose i} \d^i(a) \d^{n-i}(b).
\end{equation}

{\bf Proof of Proposition \ref{A23Dec19}}. 1. By the   definition,  the set $\{ N(E)_i\}_{i\geq 0}$ is an ascending filtration of $A$-modules on the algebra $N(E)$, and $A\subseteq N(E)_0=E^\D$. For all $i,j\geq 0$, 
$$\D^{i+j+1} (N(E)_i N(E)_j) \subseteq \sum_{s+t=i+j+1}\D^s (N(E)_i)\D^t(N(E)_j)) =0,$$
 that is $N(E)_iN(E)_j\subseteq N(E)_{i+j}$.
 
2. Statement 2 is obvious. $\Box $\\

{\bf Proof of Theorem \ref{B23Dec19}.} 1. Let $N=N_\D (E)$,   $N_i=N_\D (E)_i$ and $\ga_i = \ga\cap  N_i=\{ a\in \ga\, | \, \D^{i+1}a=0\}$ for all $i\geq 0$. Then 
$$\ga = \bigcup_{i\geq 0} \ga_i.$$

(i) $\ga_0\neq 0$ {\em is a nonzero ideal of the algebra $N_0$ (such that $\D \ga_0=0$)}: Since $\ga \neq 0$ and $\ga = \bigcup_{i\geq 0} \ga_i$, we must have $\ga_n\neq \ga_{n-1}$ for some $n\geq 0$, e.g. $n=\min \{ i\geq 0\, | \, \ga_i\neq 0\}$. Then
$$ 0\neq \D^n\ga_n\subseteq \D^nN_n\cap \ga\subseteq N_{n-n}\cap \ga = N_0\cap \ga = \ga_0.$$
(ii) $\ga'\cap N_0=\ga_0$: Notice that $\ga'\subseteq \ga$. Then 
$$ \ga_0\subseteq \ga'\cap N_0\subseteq \ga \cap N_0=\ga_0,$$
and the statement (ii) follows.

2. $\D [N_1,\ga_0]\subseteq  [\D N_1,\ga_0]+ [N_1,\D \ga_0]\subseteq [N_{1-1},\ga_0]+ [N_1,0]= [N_0,\ga_0]=0$ since the algebra $N_0$ is a commutative algebra. $\Box $\\

{\bf Proof of Theorem \ref{23Dec19}.} Since the algebra $A$ is a commutative algebra, we have that $\D \subseteq \Der_A(E)$, $A\subseteq N_\D (E)_0$,  and every ideal of the algebra $N_\D (E)$ is $\D$-stable (by the choice of $\D$). Now, the theorem follows from Theorem \ref{B23Dec19}. $\Box$

 \begin{corollary}\label{a23Dec19}
Let $A$ be a commutative subalgebra of an algebra $E$, $\D =\{ \ad_a\, | \, a\in A'\}$ where $A'$ is a non-empty subset of $A$ (eg, $A'=A$), and $R$ be a subalgebra of $N_\D (E)$ such that $A\subseteq R$. Then
\begin{enumerate}
\item The algebra $R$ is a $\D$-locally nilpotent algebra and $\{ R_i:=R\cap N_\D (E)_i\}_{i\geq 0}$ is its $\D$-order filtration. 
\item If $\ga$ is a nonzero ideal of the algebra $R$ then $\ga_0:=\ga\cap R_0=\ga\cap R^\D$ is a nonzero ideal of the algebra $R_0=R^\D\supseteq A$ such that $R\ga_0 R\cap R_0=\ga_0$. 
\item If, in addition, the algebra $R_0=R^\D$ is a commutative algebra then $[R_1, \ga_0]\subseteq \ga_0$. 
\end{enumerate}
\end{corollary}

{\it Proof.} The corollary follows from Theorem \ref{23Dec19}. $\Box$\\

{\bf Proof of Theorem \ref{C23Dec19}.} The algebra $R$ is a $\D$-locally nilpotent algebra where $\D = \{ \ad_a\, | \, a\in  A\}$ such that $R_0=R\cap \CD (A)_0= R\cap A = A$ is a commutative algebra and $\{R_i\}_{i\geq 0}$ is the $\D$-order filtration on $R$. Now, the theorem follows from Corollary \ref{a23Dec19}. $\Box$

\begin{corollary}\label{aC23Dec19}
Let $A$ be a commutative algebra, $\CD (A)$ be the algebra of differential operators on $A$, $R$ be a subalgebra  of $\CD (A)$ 
 that is generated by the algebra $A$ and a non-empty subset $\Xi$ of $ \Der_K(A)$. 
 Then $R$ is a simple algebra iff the algebra $A$ is $\Xi$-simple (i.e. $0$ and $A$ are the only $\Xi$-stable ideals of the algebra $A$). In particular, the algebra $A$ is a domain provided it is a Noetherian algebra. 
\end{corollary}

{\it Proof}. (i) {\em The algebra $R$ is not simple $\Rightarrow$ the algebra $A$ is not $\Xi$-simple}: The algebra $R$ is a subalgebra of $\CD (A )$ such that $A \subseteq R$. Suppose that $I$ is a proper ideal of $R$. Then, by Theorem \ref{C23Dec19}, $I\cap A$ is a proper $\G$-stable ideal of $A$ where $\G = R \cap \Der_K(R )\supseteq \Xi$. So, the intersection $I\cap A$ is a proper $\Xi$-stable ideal of $A$, i.e. the  algebra $A$ is not  $\Xi$-simple. 

(ii) {\em The algebra $A$ is not $\Xi$-simple $\Rightarrow$ the  algebra $R$ is not a  simple algebra}: If $J$ is a proper $\Xi$-stable ideal of $A$ then $JR$ is a proper ideal of $R$ (since $\Xi (J ) \subseteq  J$). 

Since for Noetherian algebra the minimal primes are derivation-stable, the  algebra $A$ must be a domain. So, the corollary follows. $\Box $\\

{\bf Generalized Weyl algebras, \cite{Bav-GWA-FA-91, Bav-SimGWA-1992, Bav-GWArep}.} Let $D$ be a ring, $\sigma=(\sigma_1,...,\sigma_n)$ be  an $n$-tuple of
commuting automorphisms of $D$,  $a=(a_1,...,a_n)$ be  an $n$-tuple  of elements of  the centre
$Z(D)$  of $D$ such that $\sigma_i(a_j)=a_j$ for all $i\neq j$. The {\bf generalized Weyl algebra} $A=D[X, Y; \sigma,a]$ (GWA) of rank  $n$  is  a  ring  generated
by $D$  and    $2n$ indeterminates $X_1,...,X_n, Y_1,...,Y_n$
subject to the defining relations:
$$Y_iX_i=a_i,\;\; X_iY_i=\sigma_i(a_i),\;\; X_id=\sigma_i(d)X_i,\;\; Y_id=\sigma_i^{-1}(d)Y_i\;\; (d \in D),$$
$$[X_i,X_j]=[X_i,Y_j]=[Y_i,Y_j]=0, \;\; {\rm for \; all}\;\; i\neq j,$$
where $[x, y]=xy-yx$. We say that  $a$  and $\sigma $ are the  sets  of
{\it defining } elements and automorphisms of the GWA $A$, respectively.\\

The GWA $A=\bigoplus_{\alpha\in \Z^n} Dv_\alpha$ is a $Z$-graded algebra ($Dv_\alpha \cdot Dv_\beta \subseteq Dv_{\alpha + \beta }$ for all $\alpha , \beta \in \Z^n$) where $v_\alpha = v_{\alpha_i}(1)\cdots v_{\alpha_n}(n)$ and $v_{\alpha_i}(i)=X_i^{\alpha_i}$  if $\alpha_i\geq 0$ and $v_{\alpha_i}(i)=Y_i^{-\alpha_i}$  if $\alpha_i\leq 0$.\\

The Weyl algebra $A_n$ is a generalized Weyl algebra
 $A=D[X, Y; \s ;a]$ of rank $n$ where
$D=K[H_1,...,H_n]$ is a polynomial ring   in $n$ variables with
 coefficients in $K$, $\s = (\s_1, \ldots , \s_n)$ where $\s_i(H_j)=H_j-\delta_{ij}$ and
 $a=(H_1, \ldots , H_n)$.  The map
$$A_n\ra A,\;\; X_i\mapsto  X_i,\;\; Y_i \mapsto  Y_i,\;\;  i=1,\ldots ,n,$$
is an algebra  isomorphism (notice that $Y_iX_i\mapsto H_i$).

Many quantum algebras of small Gelfand-Kirillov dimension are GWAs (eg, $U({\rm sl}_2)$, $U_q({\rm sl}_2)$, the quantum Weyl algebra, the quantum plane, the Heisenberg algebra and its quantum analogues, the quantum sphere,  and many others).\\

In case of GWAs of rank 1 we drop the lower index `1'. So, a GWA of rank 1 $A=D[x,y;\s ,a]$ is generated by the algebra $D$, $x$ and $y$ subject to the defining relations:
$$yx=a,\;\; xy=\s (a),\;\; xd=\sigma (d)x \;\; {\rm and}\;\; yd=\sigma^{-1}(d)y\;\; (d \in D).$$
The algebra $A=\bigoplus_{i\in \Z}Dv_i$ is a $\Z$-graded algebra ($Dv_iDv_j\subseteq Dv_{i+j}$ for all $i,j\in \Z$) where $v_0=1$, $v_i=x^i$ and $v_{-i}=y^i$ for $i\geq 1$.
In particular, the (first) Weyl algebra $A_i=K\langle x,\der \, | \, \der x-x\der =1\rangle\simeq K[h][x,\der; \s , a=h+1]$ is a GWA where $h=x\der$ and $\s (h) = h-1$ (since $a=\der x=x\der +1=h+1$ and  $xh=xx\der=x(\der x+[x,\der ])=(h-1)x$). \\

Let $A=D[X,Y; \s , a]$  be a GWA of rank $n$ where $D$ is  a $K$-algebra. The algebra $A$ contains two polynomial subalgebras $P_n=K[X_1, \ldots , X_n]$ and $P_n'=K[Y_1, \ldots , Y_n]$  in $n$ variables. Recall that $\s = (\s_1, \ldots , \s_n)$ where $\s_i$ are commuting automorphisms of the algebra $D$. Notice that $\s_i^{\pm 1}-1$ is a $\s_i^{\pm 1}$-derivation of the algebra $D$. A $K$-linear map $\d :D\ra D$ is called a $\s_i$-{\em derivation} if  
$$\d (ab) = \d (a) b+\s_i (a ) \d (b)\;\; {\rm  for\; all\; elements}\;\;a,b\in D.$$ 
The sets $\D =\{ \ad_{X_1}, \ldots , \ad_{X_n}\}$ (resp., $\D' =\{ \ad_{Y_1}, \ldots , \ad_{Y_n}\}$) consists of commuting $P_n$-derivations (resp., ($P_n'$-derivations) of $A$. Since $(\s_i^{-1}-1)^m=(-1)^m\s_i^{-m}(\s_i-1)^m$ for all $m\geq 1$, the map $\s_i-1$ is a locally nilpotent map on $D$ iff so is the map  $\s_i^{-1}-1$.

\begin{lemma}\label{a1Jan20}
Let $A=D[X,Y; \s , a]$  be a GWA of rank $n$ where $D$ is  a $K$-algebra, $\D =\{ \ad_{X_1}, \ldots , \ad_{X_n}\}$ and $\D' =\{ \ad_{Y_1}, \ldots , \ad_{Y_n}\}$.  
\begin{enumerate}
\item The algebra  $A$ is a $\D$-locally nilpotent algebra iff the maps $\s_1-1, \ldots , \s_n-1$ are locally nilpotent maps on $D$. 
\item The algebra  $A$ is a $\D'$-locally nilpotent algebra iff the maps $\s_1^{-1}-1, \ldots , \s_n^{-1}-1$ are locally nilpotent maps on $D$. 
\end{enumerate}
\end{lemma}

{\it Proof}. 1. $(\Rightarrow )$ The implication follows form the equality
$\ad_{X_i}^m(d)=(\s_i-1)^m(d)X_i^m$ for all $i=1, \ldots , n$ and $d\in D$ and the fact that the algebra $A$ is a $\Z$-graded algebra.  

$(\Leftarrow )$ Given an element $dv_\alpha\in Dv_\alpha$, where $d\in D$ and $\alpha =(\alpha_1, \ldots , \alpha_n)\in \Z^n$, we have to show that $\ad_{X_i}^\beta (dv_\alpha)=0$ for some element $\beta\geq 1$, by Lemma \ref{b1Jan20}. Using the $\Z$-grading of the GWA $A$, we may assume that  $\alpha_i\geq 0$. Then 
$$\ad_{X_i}^\beta (dX^\alpha)=\ad_{X_i}^\beta (d)X^\alpha=(\s_i-1)^\beta  (d)X_i^\beta X^{\alpha}  
$$
and the result follows since the map $\s_i-1$ is a locally nilpotent map on $D$.

2. Statement 2 can be proven in a similar way as statement 1. $\Box $\\

An element $r$ of a ring $R$ is called a {\bf normal element} if $rR=Rr$, i.e. $(r)=rR=Rr$ is an ideal of $R$. Given an element $s\in R$. If the set $S_s=\{ s^i \, | \, i\geq 0\}$ is a left denominator set of the ring $R$ then we denote by $R_s$ the localization $S_s^{-1}R$ of the ring $R$ at the powers of the element $s$. 

Proposition \ref{A25Dec19} is about properties of three $\D$-locally nilpotent subalgebras, $R_i$ ($i=0,1,2$),  of the  Weyl algebra $A_i$. These algebras are not simple and have very different ideal structure.

\begin{proposition}\label{A25Dec19}
Let $A_1=K\langle x,\der\rangle$ be the Weyl algebra over a field $K$ of characteristic zero and $K[x]\subset R_0\subset R_1\subset R_2\subset A_1$ be subalgebras of $A_1$ where $R_0=K\langle h=x\der , x\rangle$, $R_1=K\langle h\der , h , x\rangle$,  and  $R_2=K\langle h\der^2 , h\der , h , x\rangle$. Then 
\begin{enumerate}
\item The algebras $R_0$, $R_1$,  and $R_2$ are non-simple, $\D$-locally nilpotent algebras where $\D = \{ \ad_a\, | \, a\in K[x]\}$.
\item The algebra $R_0=K[h][x; \s ]$ is a skew polynomial ring where $\s\in \Aut_K(K[x])$ and $\s (h)=h-1$; the element $x$ of $R_0$ is a normal element; $\Spec (R_0)=\{ 0, (x), (x, p) \, | \, p\in \Irr_1 (K[h])\}$ where $\Irr_1 (K[h])$ is the set of monic irreducible polynomials of $K[h]$ (monic means that the leading coefficient of the polynomial is 1).
\item The algebra $R_1=K[h][x, y=h\der ; \s , a=h(h+1)]$ is a GWA where  $\s (h)=h-1$. The ideal $\gm_1=(y,h,x)$ is the only proper ideal of the algebra $R_1$, $R_1=K\oplus \gm_1$, $R_1/\gm_1=K$, the ideal $\gm_1$ is a maximal ideal of $R_1$ such that $\gm_1^2=\gm_1$, and $\Spec (R_1)=\{ 0, \gm_1\}$.
\item The algebra $R_2=\bigoplus_{i\geq 1} K[h]h\der^i\oplus R_0$ is a maximal subalgebra of the of the Weyl algebra $A_1$. The ideal $\gm_2=(h)=\bigoplus_{i\geq 1} K[h]h\der^i\oplus hK[h]\oplus \bigoplus_{i\geq 1} K[h]x^i$ is the only  proper ideal of the algebra $R_2$, $R_2=K\oplus \gm_2$, $R_2/\gm_2=K$, the ideal $\gm_2$ is a maximal ideal of $R_2$ such that $\gm_2^2=\gm_2$, and $\Spec (R_2)=\{ 0, \gm_2\}$.
\end{enumerate}
\end{proposition}

{\it Proof}. By the very definition, the algebras $R_0$, $R_1$ and $R_2$ are homogeneous subalgebras of the Weyl algebra $A_1$ (with respect to the $\Z$-grading of $A_1$ as a GWA). 

1. Since $A_1\simeq \CD (K[x])$ and $K[x]\subseteq R_0\subseteq R_1\subseteq R_2\subseteq A_1$, statement 1 follows.

2. It is obvious that  $R_0=K[h][x; \s ]=\bigoplus_{i\geq 0} K[h]x^i$ is a skew polynomial ring, where $\s\in \Aut_K(K[x])$ and $\s (h)=h-1$, and the element $x$ is a normal element of $R_0$ such that $R_0/(x)\simeq K[h]$. In particular, the ideal $(x)$ is a proper, prime ideal of $R_0$. The $\D$-order filtration on $R_0$ is $\{ R_{0,i}=\bigoplus_{j=0}^iK[x]h^j\}_{i\geq 0}$ since the algebra $$R_0=K[x][h; x\frac{d}{dx}]=\bigcup_{i\geq 0}R_{0,i}$$
 is an Ore extension. In particular, $R_{0,1}=K[x]\oplus K[x]h$ and for all polynomials $ p\in K[x]$, $[h, p]=x\frac{dp}{dx}$. The derivation $x\frac{d}{dx}$  is a {\em semi-simple} derivation of the polynomial algebra $K[x]=\bigoplus_{i\geq 0}Kx^i$ since $x\frac{dx^i}{dx}=ix^i$ for all $i\geq 0$. Therefore, $\{ x^iK[x]\,  | \,  i\geq 0\}$ is  the  set of $x\frac{d}{dx}$-stable ideals of the polynomial algebra $K[x]$.  The algebra $$R_{0,x}\simeq A_{1, x}$$ is a simple algebra. Hence,  if $I$ is a nonzero ideal of the algebra $R_0$ then $x^iK[x]\subseteq I$ for some $i\geq 0$. If $\gp$ is a nonzero prime ideal of $R_0$ then $\gp \supseteq (x^i)=(x)^i$ since $x$ is a normal element of $R_0$. Hence, $$(x)\subseteq \gp .$$ If $(x)\neq \gp$ then $\gp = (x,p)$ for some element $p\in \Irr_1(K[h])$ since $R_0/(x)=K[h]$. Hence $\Spec (R_0)=\{ 0, (x), (x, p) \, | \, p\in \Irr_1 (K[h])\}$ since $R_0$ is a domain.

3.  Since $yx=a$, $xy=\s (a)$,  $xd=\sigma (d)x$ and  $ yd=\sigma^{-1}(d)y$ for all $d \in D$, there is an algebra epimorphism 
$$K[h][x,y;\s , a]\ra R_1, \;\; h\mapsto h, \;\; x\mapsto x, \;\; y\mapsto h\der$$
which is an isomorphism since  $R_1=\bigoplus_{i\geq 1}K[h]y^i\oplus \bigoplus_{j\geq 0}K[h]x^j$. By \cite[Theorem 5]{Bav-UMJ-Ids-I} or \cite[Proposition 6]{Bav-UMJ-Ids-I}, the ideal $\gn_1$ is the only proper ideal of the algebra $R_1$ and $\gn_1^2=\gm_1$. Clearly, $R_1=k\oplus \gm_1$ and $R_1/\gm_1=K$, and so $\Spec (R_1)=\{ 0 , \gm_1\}$.

4. (i) $R_2=\bigoplus_{i\geq 1} K[h]h\der^i\oplus R_0$: Notice that $$A_1=\bigoplus_{i\geq 1} K[h]\der^i\oplus \bigoplus_{i\geq 0} K[h]x^i=\bigoplus_{i\geq 1} K[h]\der^i\oplus R_0$$ and $[h\der^i, h\der^j]=(i-j)h\der^{i+j}$ for all $i,j\geq 0$. Since $h\der, h\der^2\in R_2$ we have that $[h\der^2, h\der ]=h\der^3\in R_2$. Now, using induction on $i\geq 0$ and the equalities 
$$[h\der^i, h\der]=(i-1)h\der^{i+1},$$
 we see that $h\der^i\in R_2$ for all $i\geq 0$. Hence, the algebra $R_2$ contains the direct sum, say $R_2'$, from the statement (i). The direct sum $R_2'$ is a subalgebra of $A_1$ which is generated by the elements $x$ and $h\der^i$ where $i\geq 0$, i.e. $R_2'=R_2$.

(iii) {\em The algebra $R_2$ is a maximal subalgebra of the Weyl algebra} $A_1$: Suppose that $A$ be a subalgebra of $A_1$ that properly contain the algebra $R_2$. We have to show that $A=A_1$. The Weyl algebra $A_1=\bigoplus_{i\geq 1} K[h]\der^i\oplus \bigoplus_{i\geq 0} K[h]x^i$ is a direct sum of distinct eigen-spaces for the inner derivation $\ad_h$ of $A_1$ (since $[h,x^i]=ix^i$ and $[h,\der^i]=-i\der^i$ for all $i\geq 0$ and char($K$)=0). Since $h\in R_2\subseteq A$, the algebra $A$ is an $\ad_h$-stable ($[h, A]\subseteq A$). So, the algebra $A$ is a homogeneous subalgebra of the Weyl algebra $A_1$. Since 
$$A_1=\bigoplus_{i\geq 1} K[h]\der^i\oplus R_0\supseteq R_2=\bigoplus_{i\geq 1} K[h]h\der^i\oplus R_0\;\; {\rm  and}\;\; K[h]=K\oplus hK[h],$$ we must have $\der^i\in A$ for some $i\geq 1$. Then
 $$ \der= \frac{1}{i!}(-\ad_x)^{i-1}(\der^i)\in A,$$
 and so $A=A_1$ since $x,\der \in A$. 
 
 (iii) $\gm_2=(h)=\bigoplus_{i\geq 1} K[h]h\der^i\oplus hK[h]\oplus \bigoplus_{i\geq 1} K[h]x^i$: The statement (iii) follows from the statement (i) and the equalities $[h,x^i]=ix^i$ and $[h,\der^i]=-i\der^i$ for all $i\geq 0$. 
 
 By the statement (iii), $R_2=K\oplus \gm_2$ and $R_2/\gm_2=K$. 
 
 (iv) {\em The set $S_x=\{ x^i\, | \, i\geq 0\}$ is a left and right denominator set of the domains $R_0$, $R_1$, $R_2$,  and $A_1$ such that} $R_{0,x}=R_{1,x}=R_{2,x}=A_{1,x}$:
By the statement (iii), 
$$\gm_2=(h)=(x,h, h\der , \ldots , h\der^i , \ldots ) = \bigoplus_{i\geq 1} K[h]h\der^i\oplus hK[h]\oplus \bigoplus_{i\geq 1} K[h]x^i,$$
 since $[h, h\der^i]=-ih\der^i$ for all $i\geq 1$. Hence, $R_2=K\oplus \gm_2$ and $R_2/\gm_2=K$, and so $\gm_2$ is a maximal ideal of the algebra $R_2$.

 The set $S_x=\{ x^i\, | \, i\geq 0\}$ is a left and right Ore set of the domains $R_0\subseteq R_1\subseteq R_2\subseteq A_1$ (use the $\Z$-gradings of the algebras). Since $\der =x^{-1} x\der\in R_{0,x}$, we see that $R_{0,x}=A_{1,x}$. Then the inclusions $R_0\subseteq R_1\subseteq R_2\subseteq A_1$ yield the equalities $R_{0,x}=R_{1,x}=R_{2,x}=A_{1,x}$. 

(v) {\em The ideal $\gm_2$ is the only proper ideal} of $R_2$: Let  $I$ be a proper ideal $I$ of $R_2$ we have to show that $I=\gm_2$. By the statement (iv), $R_1\subset R_2\subset R_{1,x}=R_{2,x}$, and so the algebra $R_1$ is an essential left $R_1$-submodule of $R_2$. Hence $I\cap R_1=\gm_1$ is the only proper ideal of the algebra $R_1$, by statement 3. Since 
$h\der \in \gm_1\subseteq \gm_2$ and $$I\ni [h\der^i, h\der ] =(i-1) h\der^{i+1}\;\; {\rm for \; all}\;\; i\geq 2,$$
 we have that  $I\supseteq (x, h, h\der, \ldots , h\der^i , \ldots ) = \gm_2$,  i.e. $I=\gm_2$, by the maximality of the ideal $\gm_2$. 
 
 The algebra $R_2$ is a domain, hence $\gm_2^2=\gm_2$, by the statement (v). Now, $\Spec (R_2) = \{ 0, \gm_2\}$.  $\Box $\\
 
 {\bf Localizations and the algebras $N_\D (E)$.} An element of a ring is called a {\em regular} element if it is not a zero divisor.  Every regular left Ore set of a ring $R$  is a regular left denominator set, and vice versa. The set of all regular left Ore sets of $R$  is denoted by $\Den_l(R, 0)$. Proposition \ref{A2Jan20} shows that the algebra $N_\D (E)$ is well-behaved under localizations at regular left Ore sets that are contained in the zero component $N_\D (E)_0$ of the $\D$-order filtration.

\begin{proposition}\label{A2Jan20}
 Let $A$ be a subalgebra of an algebra $E$, $\D\subseteq \Der_A(E)$, and $S\in \Den_l(E,0)$ with $S\subseteq N_\D (E)_0$. Then 
\begin{enumerate}
\item $A\subseteq E\subseteq S^{-1}E$ and $ \D\subseteq \Der_A(S^{-1}E)$. 
\item $S\in \Den_l(N_\D (E),0)$.
\item $N_\D (S^{-1}E)\simeq S^{-1}N_\D (E)$.
\item For all integers $i\geq 0$, $N_\D (S^{-1}E)_i\simeq S^{-1}N_\D (E)_i$.
\end{enumerate}
\end{proposition}

{\it Proof}. Let $N=N_\D (E)$ and $N_i=N_\D (E)_i$ for $i\geq 0$.

1. Statement 1 is obvious (for all elements $\d\in \D$, $s\in S$ and $e\in E$, $\d (s^{-1}e)=s^{-1}\d (e)$ since $\d (s)=0$). 

2. Clearly, $S\subseteq N_0\subseteq N$. We have to show that the set $S$ is a left Ore set of $N$. Given elements $s\in S$ and $ n\in N$, i.e. $\D^in=0$ for some $i\geq 1$. Then $ns^{-1} = t^{-1}e$  for some elements $ t\in S$ and $ e\in E$. Then
$$ 0=\D^i(n)s^{-1}=\D^i(ns^{-1})=\D^i ( t^{-1}e)=t^{-1}\D^i(e),$$
and so $\D^i(e)=0$, that is $e\in N$. Therefore, $tn=es$. This means that the set $S$ is a left Ore set in $N$.

4. Clearly, $S^{-1}N\subseteq N_\D (S^{-1}E)$. Given an element $s^{-1}n\in N_\D (S^{-1}E)$. Then $0=\D^i(s^{-1}n) = s^{-1}\D^i(n)$ iff $\D^i(n)=0$, and statement 4 follows.

3. Statement 3 follows from statement 4.  $\Box $\\

Let $I$ be an ideal of a ring $E$. We denote by $\Den_l(E, I)$ the set of left denominator sets of $E$ with $I=\ass_E(S):=\{ e\in E\, | \, se=0$ for some element $s\in \}$. In case $S\in \Den_l(E, I)$ and $I\neq 0$, we have to impose an additional condition that $|\D |<\infty$ (the set $\D$ is a finite set) in order to have similar results as in Proposition \ref{A2Jan20}, see Proposition \ref{B2Jan20}.

\begin{proposition}\label{B2Jan20}
 Let $A$ be a subalgebra of an algebra $E$, $\D\subseteq \Der_A(E)$, and $S\in \Den_l(E,I)$ with $S\subseteq N_\D (E)_0$, $\bE =E/I$, $\bA = A/I'$ where $I'=A\cap I$,   $\overline{\D}=\{ \overline{\d}\, | \, \d \in \D\}\subseteq \Der_{\bA}(\bE )$,  and  $\overline{\d}(e+I)=\d (e)+I$ for all elements $e\in E$ (see statement 1). Then 
\begin{enumerate}
\item The ideal $I$ is $\D$-stable ($\D I\subseteq I$).
\item $S\in \Den_l(N_\D (E), N_\D (E)\cap I)$ provided $|\D | <\infty$.
\item $N_\D (S^{-1}E)\simeq S^{-1}N_\D (E)$  provided $|\D | <\infty$.
\item For all integers $i\geq 0$, $N_\D (S^{-1}E)_i\simeq S^{-1}N_\D (E)_i$  provided $|\D | <\infty$.
\item $\bS =\{ s+I \, | \, s\in S\}$,    $\bS^{-1}\bE\simeq  S^{-1}E$,  $\bS^{-1} N_{\overline{\D}}(\bE )\simeq N_{\overline{\D}}(\bS^{-1}\bE )\simeq N_\D (S^{-1}E)$,  and  $\bS^{-1} N_{\overline{\D}}(\bE )_i\simeq N_{\overline{\D}}(\bS^{-1}\bE )_i\simeq N_\D (S^{-1}E)_i$ for all $i\geq  0$.
\end{enumerate}
\end{proposition}

{\it Proof}. Let $N=N_\D (E)$ and $N_i=N_\D (E)_i$ for $i\geq 0$.

1. Given elements $\d \in \D$ and $a\in I$. Then $sa=0$ for some element $s\in S$, and so $ 0=\d(sa)=s\d (a)$. This implies that $\d (a)\in I$, and statement 1 follows. 

2. (i) $S$ {\em is a left Ore set of} $N$: Given elements $s\in S$ and $n\in N$, we have to show that $s_1n=n_1s$ for some elements $s_1\in S$ and $n_1\in N$. Since $S$ is a left Ore set of $E$, $tn=es$ for some elements $t\in S$ and $e\in E$. Since $n\in N$, $\D^in=0$ for some $i\geq 1$. Then 
$$0=t\D^in=\D^i(tn)=\D^i(es)=\D^i(e)s,$$
and so $\D^i(e)\subseteq I$ (since $S\in \Den_l(E,I)$). The set $\D$ is a finite set hence so is the set $\D^i(e)$. We can fix an element $s'\in S$ such that $0=s'\D^i(e)=\D^i(s'e)$, i.e. $n_1:=s'e\in N_{i-1}$. Now, it suffices to take $s_1=s't$ since $$s_1n=s'tn=s'es=n_1s.$$

(ii) $S\in \Den_l(N, N\cap I)$: Since $\ass_E(S)=I$, we have that $\ass_N(S)=N\cap I$. If $ns=0$ for some elements  $n\in N$ and $s\in S$. Then $n\in N\cap I$, and the statement (ii) follows from the statement (i).

4. (i) $S^{-1}N_i\subseteq N(S^{-1}E)_i$ {\em for all } $i\geq 0$: Given elements $s\in S$ and $n\in N_i$. Then $\D^{i+1}n=0$ and $$\D^{i+1}(s^{-1}n)=s^{-1}\D^{i+1}n=0,$$ and so $s^{-1}n\in  N(S^{-1}E)_i$.

(ii) $N(S^{-1}E)_i\subseteq S^{-1}N_i $ {\em for all } $i\geq 0$: Given an element $t^{-1}e\in N (S^{-1}E)_i$ where $t\in S$ and $e\in E$. Then $$0=\D^{i+1}(t^{-1}e) = t^{-1}\D^{i+1}(e),$$ and so the set $\D^{i+1}(e)$ is a finite subset of $I$ (since $|\D |<\infty$). Hence, there exists an element $t_1\in S$ such that $0=t_1\D^{i+1}(e)=\D^{i+1}(t_1e)$, i.e. $t_1e\in N_i$. Now, $t^{-1}e=(t_1t)^{-1}t_1e\in S^{-1}N_i$, and the statement (ii) follows.

3. Statement 3 follows from statement 4.  

5. Statement 5 follows from Proposition \ref{A2Jan20}. $\Box $\\

{\bf Monoids that are generated by $\ad$-locally nilpotent elements are denominator sets.} Let $R$ be a ring and $s,r\in R$.  Then 
\begin{equation}\label{smr=mi}
s^mr=\sum_{i=0}^m {m\choose i} \ad_s^i(r) s^{m-i}\;\; {\rm for\; all}\;\; m\geq 1,
\end{equation}
\begin{equation}\label{smr=mi1}
rs^m=\sum_{i=0}^m {m\choose i}  s^{m-i} (-\ad_s)^i(r)\;\; {\rm for\; all}\;\; m\geq 1.
\end{equation}
Suppose that $rs=0$ (resp., $sr=0$) then by (\ref{smr=mi}) ( resp., (\ref{smr=mi1})) for all $n\geq 1$, 
\begin{equation}\label{snr=ads}
s^nr=\ad_s^n(r) \;\;\; ({\rm resp.,}\; rs^n=(-\ad_s)^n(r)).
\end{equation}

Let $R$ be a ring and $S\in \Den_l(R, \ga )$. A (left) ideal $I$ of $R$ is called an $S$-{\em saturated ideal} if the inclusion $sr\in I$ (where $s\in S$ and $r\in R$) implies the inclusion $r\in I$, i.e. $$\tor_S(R/I):=\{ a\in R/I\, | \, sa=0\;\; {\rm for\; some}\;\; s\in S \}=0,$$ the $R$-module $R/I$ is $S$-torsionfree. In general, if $I$ is an ideal of $R$ the localization $S^{-1}I$ (which is a left ideal of the ring $S^{-1}R$) is not an ideal of $S^{-1}R$. Proposition \ref{C2Jan20}.(2) gives a class of denominator sets $S$ of an arbitrary ring $R$  such that  $S^{-1}I$ is always an ideal of $S^{-1}R$. We denote by $\CI (R)$ and $\CI (R, S-{\rm sat.})$ the sets of ideals and $S$-saturated ideals of the ring $R$, respectively. 

\begin{proposition}\label{C2Jan20}
Let $R$ be a ring and $S$ be a multiplicative subset of $R$. Suppose that the monoid $S$ is generated by a set of $\ad$-locally nilpotent elements, say $S=\langle s_\l\, | \, \l \in \L\rangle$ (the inner derivations $\{ \ad_{s_\l}\, | \, \l \in \L\}$ of $R$ are locally nilpotent). Then 
\begin{enumerate}
\item $S\in \Den (R, \ga  )$.
\item If $I$ is an ideal of the ring $R$ then $S^{-1}I=IS^{-1}$ is an ideal of the ring $S^{-1}R\simeq RS^{-1}$.
\item The map $ \CI (R, S-{\rm sat.})\ra \CI (S^{-1}R)$, $I\mapsto S^{-1}I$ is a bijection with the inverse $J\mapsto \s^{-1}(J)$ where $\s : R\ra S^{-1}R$, $r\mapsto \frac{r}{1}$.
\end{enumerate}
\end{proposition}

{\it Proof}. 1. (i) $S$ {\em is an Ore set of} $R$: To prove the statement (i) it suffices to show that the left (resp., right) Ore condition holds for the generators $\{ s_\l\}$ of the monoid $S$. Since the maps $\ad_{s_\l}$ are locally nilpotent, this follows from Eq. (\ref{smr=mi}) (resp.,  Eq. (\ref{smr=mi1})).

(ii) $S\in \Den (R, \ga )$: If $rs_\l =0$ (resp., $s_\l r=0$) for some $\l \in \L$ and $r\in R$ then, by Eq. (\ref{snr=ads}), $s_\l^{n(\l )}r=0$ (resp., $rs_\l^{n(\l )}=0$) for some natural number $n(\l )$. Using this fact we see that if $rs_\l \cdots s_\mu =0$ (resp., $s_\l \cdots s_\mu r=0$) then $s_\l^{n(\l )}\cdots s_\mu^{n(\mu )}r=0$ (resp., $r s_\l^{n(\l )}\cdots s_\mu^{n(\mu )}=0$ and the statement (ii) follows (recall that every element $s\in S$ is a product $s=s_\l \cdots s_\mu$),

2. Let $I$ be an ideal of the ring $R$. By statement 1, $S^{-1}I$ (resp., $IS^{-1}$) is a left (resp., right) ideal of the ring $S^{-1}R$ (resp., $RS^{-1}$). Since $S\in \Den (R)$, $S^{-1}R\subseteq RS^{-1}$. The inclusion $S^{-1}I\subseteq IS^{-1}$ (resp., $IS^{-1}\subseteq S^{-1}I$) follows from the equality: For all elements $s_\l $ and $r\in R$,
$$ s_\l^{-1}r=s_\l^{-1}(rs_\l^n)s_\l^{-n}=s_\l^{-1}\bigg( \sum_{i=0}^n {n\choose i}s_\l^i (-\ad_{s_\l})^{n-i}(r)\bigg) s_\l^{-n}\in IS^{-1}\;\; {\rm for\; all}\;\; n\gg 0,$$
 $$ ({\rm resp.,}\;\; rs_\l^{-1}=s_\l^{-n}(s_\l^nr)s_\l^{-1}=s_\l^{-n}\bigg( \sum_{i=0}^n {n\choose i} \ad_{s_\l}^i(r)s_\l^{n-i}\bigg) s_\l^{-1}\in S^{-1}I\;\; {\rm for\; all}\;\; n\gg 0).$$
3. By statement 2, the map $ \CI (R, S-{\rm sat.})\ra \CI (S^{-1}R)$, $I\mapsto S^{-1}I$ is well-defined. By the very definition, the map $ \CI (S^{-1}R)\ra \CI (R, S-{\rm sat.})$, $J\mapsto \s^{-1}(J)$ is well-defined. Since $S^{-1}\s^{-1}(J)=J$ and $\s^{-1}(S^{-1}I)=I$, statement 3 follows. $\Box $


\section{Simplicity criteria for subalgebras of $\CD (\CA )$ that contain $\CA$}\label{SIMCRIT}

The aim of the section is to prove  Theorem \ref{1Jan20}, Theorem \ref{27Dec19},   Theorem \ref{28Dec19},  and Theorem \ref{8Jan20}. Each commutative algebra $A$ is a left $\CD (A)$-module, its submodule structure is described in Proposition \ref{A3Jan20}. 
Theorem  \ref{4Jan20} gives the canonical form for each differential operator on arbitrary commutative algebra. 

The following notation will remain  fixed throughout the section (if
it is not stated otherwise): $K$ is a field of arbitrary  characteristic 
(not necessarily algebraically closed), 
 $P_n=K[x_1, \ldots , x_n]$ is a polynomial algebra over
$K$, $\der_1:=\frac{\der}{\der x_1}, \ldots ,
\der_n:=\frac{\der}{\der x_n}\in \Der_K(P_n)$, $I:=\sum_{i=1}^m
P_nf_i$ is a {\bf prime} but {\bf not} a maximal ideal of the
polynomial algebra $P_n$ with a set of generators $f_1, \ldots ,
f_m$, the algebra $A:=P_n/I$ which is a domain with the field of
fractions $Q:={\rm Frac}(A)$, the epimorphism $\pi :P_n\ra A$,
$p\mapsto \bp :=p+I$, to make notation simpler we sometime write
$x_i$ for $\overline{x}_i$ (if it does not lead to confusion), the
{\bf Jacobi} $m\times n$ matrices
 $$J=\bigg(\frac{\der f_i}{\der x_j}\bigg)\in M_{m,n}(P_n)$$ and
 $\bJ =\bigg(\overline{\frac{\der f_i}{\der
x_j}}\bigg)\in M_{m,n}(A)\subseteq M_{m,n}(Q)$,  $r:={\rm rk}_Q(\bJ )$
is the {\bf rank} of the Jacobi matrix $\bJ$ over the field $Q$,
$\ga_r$ is the {\bf Jacobian ideal} of the algebra $A$ which is
(by definition) generated by all the $r\times r$ minors of the
Jacobi matrix $\bJ$.  

 For $\i =(i_1, \ldots , i_r)$ such that $1\leq
i_1<\cdots <i_r\leq m$ and $\j =(j_1, \ldots , j_r)$ such that
$1\leq j_1<\cdots <j_r\leq n$, $\D
 (\i , \j )$ denotes the corresponding minor of the Jacobi matrix
$\bJ =(\bJ_{ij})$, that is $\det (\bJ_{i_\nu, j_\mu})$, $\nu , \mu
=1, \ldots, r$,  and the element $\i$ (resp., $\j $) is called {\bf
non-singular} if $\D (\i , \j')\neq 0$ (resp., $\D (\i', \j )\neq
0$) for some $\j'$ (resp., $\i'$). We denote by $\II_r$ (resp.,
$\JJ_r$) the set of all the non-singular $r$-tuples $\i$ (resp.,
$\j $).

Since $r$ is the rank of the Jacobi matrix $\bJ$, it is easy to
show that $\D (\i , \j )\neq 0$ iff $\i\in \II_r$ and $\j\in
\JJ_r$, \cite[Lemma 2.1]{gendifreg}. \\

 A localization of an {\em affine} algebra is
called an algebra of {\bf essentially finite type}. Let $\CA :=S^{-1} A$ be a localization of the algebra $A=P_n/I$ at a
multiplicatively closed subset $S$ of $A$. Suppose that $K$ is a perfect field. Then the algebra $\CA$ is {\em regular} iff $\ga_r=\CA$ where $\ga_r$ is the Jacobian ideal of $\CA$  (the
{\bf Jacobian criterion of regularity}, \cite[Theorem 16.19]{Eisenbook}).
For any regular algebra $\CA$ over a perfect field,  explicit sets of generators and defining relations for the algebra $\CD (\CA )$ are given in \cite{gendifreg} (char($K$)=0) and  \cite{gendifregcharp} (char($K)>0$).\\

Let $R$ be a commutative $K$-algebra. The ring of ($K$-linear)
{\bf differential operators} $\CD (R)$ on $R$ is defined as a
union of $R$-modules  $\CD (R)=\bigcup_{i=0}^\infty \,\CD (R)_i$
where 
$$ \CD (R)_i=\{ u\in \End_K(R)\, |\, [r,u]:=ru-ur\in \CD (R)_{i-1} \; {\rm for \; all \; }\; r\in R\},\;\; i\geq 0, \;\; \CD(R)_{-1}:=0.$$
In particular, $\CD (R)_0=\End_R(R)\simeq R$, $(x\mapsto bx)\lra b$.
 The set of $R$-bimodules $\{ \CD (R)_i\}_{i\geq 0}$ is the {\bf order filtration} for
the algebra $\CD (R)$:
$$\CD(R)_0\subseteq   \CD (R)_1\subseteq \cdots \subseteq
\CD (R)_i\subseteq \cdots\;\; {\rm and}\;\; \CD (R)_i\CD
(R)_j\subseteq \CD (R)_{i+j} \;\; {\rm for\; all}\;\; i,j\geq 0.$$

The subalgebra $\D (R)$ of $\CD (R)$ which is  generated by $R\equiv
\End_R(R)$ and the set ${\rm Der}_K (R)$ of all $K$-derivations of
$R$ is called the {\bf derivation ring} of $R$.

Suppose that $R$ is a  regular affine  domain of Krull dimension
$n\geq 1 $ and char($K$)=0. In geometric terms, $R$ is the coordinate ring $\OO
(X)$ of a smooth irreducible  affine algebraic variety $X$ of
dimension $n$. Then
\begin{itemize}
\item ${\rm Der}_K (R)$ {\em is a finitely generated projective}
$R$-{\em module of rank} $n$, \item  $\CD (R)=\Delta (R) $, \item
$\CD (R)$ {\em is a simple (left and right) Noetherian domain of
Gelfand-Kirillov dimension}  $\GK \, \CD (R)=2n$ ($n=\GK (R)=\Kdim
(R))$.
\end{itemize}

For the proofs of the statements above the reader is referred to
\cite{MR}, Chapter 15.
 So, the domain $\CD (R)$ is a simple finitely generated infinite dimensional Noetherian algebra
(\cite{MR}, Chapter 15).  

If char($K)>0$ then  $\CD (R)\neq\Delta (R) $ and the algebra $\CD (R$ is not finitely generated and neither left nor right Noetherian but analogues of the results above hold but the Gelfand-Kirillov dimension has to replaced by a new dimension introduced in \cite{BernIneqcharp}.

\begin{lemma}\label{b28Dec19}
Let $\CA$ be a commutative algebra of essentially finite type, $\CD (\CA )$ be the algebra of differential operators on $\CA$, $R$ be a subalgebra of $\CD (\CA )$ that contains $\CA$. Then, for every $i\geq 0$, the left and right $\CA$-module $R_i=R\cap \CD (\CA )_i$ is finitely generated and Noetherian.
\end{lemma}

{\it Proof}. For each $i\geq 0$, the left and right $\CA$-module $\CD (\CA )_i$ is finitely generated, hence Noetherian since the algebra $\CA$ is Noetherian. Since $R_i$ is a left and right $\CA$-submodule of $\CD (\CA )_i$, it  is also finitely generated and Noetherian.  $\Box $\\

The next obvious lemma is a criterion for a subalgebra of $\CD (\CA )$ that contains $\CA$ to be  an essential left or right $\CA$-submodule of $\CD (\CA )$.

\begin{lemma}\label{a27Dec19}
Let $\CA$ be a commutative domain of essentially finite type over a field of characteristic zero and $Q$ be its field of fractions, $R$ be a subalgebra of  $\CD (\CA )$ that contains $\CA$, and  $R_i=R\cap \CD (\CA )_i$ where $i\geq 0$. Then the following statements are equivalent:
\begin{enumerate}
\item $Q\t_AR=Q\t_A\CD (\CA )$.
\item $R\t_AQ=\CD (\CA ) \t_AQ$.
\item $\dim_Q(Q\t_AR)=\dim_Q(Q\t_A\CD (\CA ))$.
\item $\dim_Q(R\t_AQ)=\dim_Q(\CD (\CA ) \t_AQ)$.
\item $\dim_Q(QD_R)=\dim_Q(Q\Der_K(\CA ))$ where $D_R:=R\cap \Der_K(\CA )$.
\item $D_R$ is an essential left $\CA$-submodule of $\Der_K(\CA )$. 
\item $R_1$ is an essential left $\CA$-submodule of $\CD (\CA )_1$.
\item $R_1$ is an essential right $\CA$-submodule of $\CD (\CA )_1$.
\item There is an natural number $i\geq 1$ such that $R_i$ is an essential left $\CA$-submodule of $\CD (\CA )_i$.
\item There is an natural number $i\geq 1$ such that $R_i$ is an essential right $\CA$-submodule of $\CD (\CA )_i$.
\end{enumerate}
\end{lemma}

{\it Proof}. Straightforward. $\Box $

\begin{lemma}\label{a28Dec19}
Let $\CA$ be a commutative domain of essentially finite type over a field of arbitrary  characteristic, $R$ be a subalgebra of  $\CD (\CA )$ that contains $\CA$, and $S$ be a multiplicative subset of $\CA \backslash \{ 0\}$.
 Then $S\in {\rm Den} (R, 0)$, $S^{-1}\CA \subseteq S^{-1}R\subseteq S^{-1}\CD (\CA )\simeq \CD (S^{-1}\CA )$. If, in addition, $ S^{-1}R= S^{-1}\CD (\CA )$ then the algebra $R$ is an essential left and right $R$-submodule of the algebra $\CD (\CA )$. 
\end{lemma}

{\it Proof.} The lemma follows at once from Proposition \ref{C2Jan20} or    Eq. (\ref{smr=mi}) and the fact that $\CD (\CA )$ is a domain. $\Box$\\

{\bf Proof of Theorem \ref{27Dec19}.} Since the $\CA$-submodule $R_i$ of $\CD (\CA )_i$ is an essential submodule and the $\CA$-module $\CD (\CA )_i$ is finitely generated, $\gb_i\neq0 $. Then, by Eq.  (\ref{smr=mi1}),  $\gc_i\neq 0$.

 $(4\Leftrightarrow 5)$ Theorem \ref{CX9Jul19}.

$(1\Rightarrow 5)$ (i) {\em The algebra $\CD (\CA )$ is simple}: Otherwise, for each proper ideal $I$ of the algebra $\CD (\CA )$ the intersection $R\cap I$ is a proper ideal of the algebra $R$ since the algebra $R$ is an essential $R$-submodule of $\CD (\CA )$, this contradicts the simplicity of the algebra $R$.

(ii) $R\gb_1^2R=R$  {\em and } $R\gb_1\cdots \gb_{i-1}\gb_i^2R=R$ {\em for} $i\geq 2$: 
 Since the ideals $\gb_i$ ($i\geq 1$) of the domain $\CA$ are non-zero, so are their products $\gb_1^2$   and  $\gb_1\cdots \gb_{i-1}\gb_i^2$, and the statement (ii) follows from the simplicity of the algebra $R$.

$(5\Rightarrow  1)$ By Theorem \ref{C23Dec19}, it suffices to show that for every  nonzero ideal $\ga$ of the algebra $\CA$ such that $[R_1, \ga ] \subseteq \ga$, the ideal $R\ga R$ is equal to $R$. Since the algebra $\CD (\CA )$ is a simple algebra, $1\in \CD (\CA ) \ga \CD (\CA )$, i.e. $1\in \CD (\CA )_i \ga \CD (\CA )_i$ for some $i\geq 0$. If $i=0$ then $1\in \CA\ga \CA =\ga$ and there is nothing to prove.

So, we assume that $i\geq 1$. If $i=1$ then 
$$ \gb_1^2= \gb_1\cdot 1 \cdot \gb_1\subseteq \gb_1\cdot \CD (\CA )_1 \ga \CD (\CA )_1\cdot \gb_1\subseteq R_1\ga \Big([\CD (\CA )_1, \gb_1]+\gb_1\CD (\CA )_1\Big) \subseteq R_1\ga ( A+R_1)=R_1\ga R_1,$$
and so $R=R\gb_1^2R\subseteq R\ga R\subseteq R$, i.e. $R\ga R = R$. If $i\geq 2$ then 
\begin{eqnarray*}
\gb_i^2\gb_{i-1}\cdots \gb_1&=& \gb_i\cdot 1 \cdot \gb_i\gb_{i-1}\cdots \gb_1\subseteq 
\gb_i\cdot  \CD (\CA )_i \ga \CD (\CA )_i\cdot \gb_i\gb_{i-1}\cdots \gb_1\\
&\subseteq & R_i \ga \Big([\CD (\CA )_i, \gb_i]+\gb_i\CD (\CA )_i\Big)\gb_{i-1}\cdots \gb_1
\\
&\subseteq & R_i\ga R_i+R_i\ga \Big(\CD (\CA )_{i-1}+R_i\Big)\gb_{i-1}\cdots \gb_1\\
& \subseteq & R_i\ga R_i+R_i\ga \Big([\CD (\CA )_{i-1}, \gb_{i-1}]+\gb_{i-1}\CD (\CA )_{i-1}\Big)\gb_{i-2}\cdots \gb_1\\
&\subseteq & R_i\ga R_i+R_i\ga (\CD (\CA )_{i-2}+R_{i-1})\gb_{i-2}\cdots \gb_1\\
&\subseteq &\cdots \subseteq R_i\ga R_i+R_i\ga (\CD (\CA )_0+R_0) \\
 &\subseteq & R_i\ga R_i. 
\end{eqnarray*}
Hence, $R=R(\gb_i^2\gb_{i-1}\cdots \gb_1)R\subseteq R\ga R\subseteq R$, i.e. $R=R\ga R$, as required. 

$(6\Leftrightarrow 7)$ Theorem \ref{CX9Jul19}.

$(5\Rightarrow  7)$ (resp., $(7\Rightarrow  5)$) Repeat the proofs of the implication 
 $(1\Rightarrow  5)$ (resp., $(5\Rightarrow 1)$)  replacing the ideals $\gb_i$ by $\gc_i$ and using right modules instead the left ones.
 
$(2\Leftrightarrow 3)$ Theorem \ref{CX9Jul19}.

$(1\Rightarrow  3)$ (i) {\em The algebra $\CD (\CA )$ is simple}: See the proof of the statement (i) in the proof of the implication $(1\Rightarrow  5)$.

(ii) $R\gb_i\gc_iR=R$ {\em for all} $i\geq 1$:  Since the ideals $\gb_i$ and $\gc_i$  ($i\geq 1$) of the domain $\CA$ are non-zero, so are their products $\gb_i\gc_i$, and the statement (ii) follows from the simplicity of the algebra $R$. 

$(3\Rightarrow  1)$  By Theorem \ref{C23Dec19}, we have to show that for every  nonzero ideal $\ga$ of the algebra $\CA$ such that $[R_1, \ga ] \subseteq \ga$, the ideal $R\ga R$ is equal to $R$. Since the algebra $\CD (\CA )$ is a simple algebra, $1\in \CD (\CA ) \ga \CD (\CA )$, i.e. $1\in \CD (\CA )_i \ga \CD (\CA )_i$ for some $i\geq 0$. Then
$$\gb_i\gc_i=\gb_i\cdot 1 \cdot\gc_i\subseteq \gb_i \CD (\CA )_i \ga \CD (\CA )_i\gc_i\subseteq R_i\ga R_i,$$
and so $R=R\gb_i\gc_iR\subseteq R\ga R \subseteq R$, i.e. $R\ga R =R$. $\Box $\\

{\bf Proof of Theorem \ref{28Dec19}.} Notice that the algebra $R$ is an essential left and right $\CA$-submodule of $\CD (\CA )$ and the theorem follows from Theorem \ref{27Dec19}. Let us give more details.

The equivalences $(2\Leftrightarrow 3)$, $(4\Leftrightarrow 5)$ and $(6\Leftrightarrow 7)$ follow from  Theorem \ref{CX9Jul19}.

$(1\Rightarrow 3, 1 \Rightarrow 5, 1 \Rightarrow 7)$ The implications follow from the fact that the algebra $\CD (\CA )$ is a simple algebra and that $s_it_i\neq 0$ for all $i\geq 1$ (since the algebra $\CA$ is a domain).

$(3\Rightarrow 1, 5 \Rightarrow 1, 7 \Rightarrow 1)$ Since $s_i\in \gb_i$ and $t_i\in \gc_i$ for all $i\geq 1$, the implications follow from Theorem \ref{27Dec19}. $\Box$\\

{\bf Proof of Theorem \ref{1Jan20}}. $(1\Rightarrow 2)$ If the algebra $\CD (\CA )$ is a simple algebra then so is the algebra $\CD (\CA )_\gm$ for all maximal ideals $\gm$ of the algebra $\CA$ that contain the Jacobian ideal $\ga_r$, by Proposition \ref{C2Jan20}.(2).

$(2\Rightarrow 1)$ Suppose that $I$ is a proper ideal of the algebra $\CD (\CA )$, we seek a contradiction. Then there is a maximal ideal $\gm$ of the algebra $\CA$ such that $I_\gm$ is a proper ideal of the algebra $\CD (\CA )_\gm$, by Proposition \ref{C2Jan20}.(2) and since the algebra $\CD (\CA )$ is a domain. Since $\CD (\CA)_\gm \simeq \CD (\CA_\gm )$, we must have that $\ga_r\subseteq \gm$ (by Theorem \ref{CX9Jul19}), a contradiction. $\Box $ \\

{\bf Proof of Theorem \ref{8Jan20}}. $(1\Rightarrow 2,4)$ If the algebra $R$ is a simple algebra then,  by Proposition \ref{C2Jan20}.(2),  so is the algebra $R_\gm$ for all maximal ideals $\gm$ of the algebra $\CA$. Since the algebra $R$ is an essential left $R$-submodule of the algebra $\CD (\CA )$ the algebra $\CD (\CA )$ must be simple (since $R$ is simple). 

$(2\Rightarrow 1)$ Suppose that $R$ is not a simple algebra, we seek a contradiction. By Theorem \ref{27Dec19}.(5), one of the ideals, say $I$, in the set  $\{ R\gb_1^2R, R\gb_1\cdots \gb_{i-1}\gb_i^2R\, | \, i\geq 2\}$ is not equal to $R$. Then the ideal $\ga = I\cap \CA$ of the algebra $\CA$ is a proper ideal that contains either the ideal $\gb_1^2$ or $\gb_1\cdots \gb_{i-1}\gb_i^2$. Then there is a maximal ideal $\gm$ of $\CA$ that contains $\ga$ and such that the ideal $I_\gm$ is a proper ideal of the algebra $R_\gm $. Clearly, the ideal $\gm$ contains  one of the ideals $\gb_i$, a contradiction.

$(2\Leftrightarrow 3)$, $(4\Leftrightarrow 5)$ These implications follow from Theorem \ref{1Jan20}.

$(4\Rightarrow 1)$ Suppose that $R$ is not a simple algebra, we seek a contradiction. By Theorem \ref{27Dec19}.(5), one of the ideals, say $I$, in the set  $\{ R\gc_1^2R, R\gc_1\cdots \gc_{i-1}\gc_i^2R\, | \, i\geq 2\}$ is not equal to $R$. Then the ideal $\ga = I\cap \CA$ of the algebra $\CA$ is a proper ideal that contains either the ideal $\gc_1^2$ or $\gc_1\cdots \gc_{i-1}\gc_i^2$. Then there is a maximal ideal $\gm$ of $\CA$ that contains $\ga$ and such that the ideal $I_\gm$ is a proper ideal of the algebra $R_\gm $. Clearly, the ideal $\gm$ contains  one of the ideals $\gc_i$, a contradiction.  $\Box $\\

{\bf The $\CD (A)$-module structure of the algebra $A$ and its simplicity criterion.} Let $A$ be an {\em arbitrary}  algebra and $\CD (A)$ be the algebra of differential operators on $A$. By the  definition of the algebra $\CD (A)$, the algebra $A$ is a faithful left $\CD (A)$-module (since $\CD (A)\subseteq \End_K(A)$). The action of elements $\d\in \CD (A)$ on the elements $a\in A$ is denoted either by $\d (a)$ or $\d *a$. Since $A\subseteq \CD (A)$, 
$$A=\CD (A)*1\simeq \CD (A)/\CD (A)_{[0]}\;\; {\rm where}\;\; \CD (A)_{[0]}:=\{ \d \in \CD (A) \, | \, \d*1=0\}$$
is the annihilator of the element 1 of the $\CD (A)$-module $A$. By the  definition, $\CD (A)_{[0]}$ is a left ideal of the algebra $\CD (A)$ such that 
\begin{equation}\label{DA=ADA0}
\CD (A)=A\oplus \CD (A)_{[0]}
\end{equation}
is a direct sum of left $A$-modules. Clearly, $\Der_K(A)\subseteq \CD (A)_{[0]}$. Notice that 
\begin{equation}\label{DD0D}
\CD (A)\CD (A)_{[0]}\CD (A)=\CD (A)*A+\CD (A)_{[0]}
\end{equation}
since
\begin{eqnarray*}
\CD (A)\CD (A)_{[0]}\CD (A) &=& \CD (A)_{[0]}\CD (A)=\CD (A)_{[0]}(A+\CD (A)_{[0]}) \\
 &=& \CD (A)*A+\CD (A)_{[0]}.
\end{eqnarray*}
 We denote by $\Sub_{\CD (A)}(A)$ the set of all left $\CD (A)$-submodules of the $\CD (A)$-module $A$. Let $\CI (A, \CD (A)-{\rm st.})$ (resp., $\CI (A, \CD (A)-{\rm st.}, \CD (A)*A)$) be the set of all $\CD (A)$-stable ideals  of $A$ (resp., that contain the ideal $\CD (A)*A$). By (\ref{DA=ADA0}), an ideal $\ga$ of $A$ is $\CD (A)$-stable iff it is $\CD (A)_{[0]}$-stable ($\CD (A)_{[0]}*\ga \subseteq \ga $). Let $\CI (\CD (A), \CD (A)_{[0]})$ be the set of ideals of the algebra $\CD (A)$ that contain $\CD (A)_{[0]}$. Proposition \ref{A3Jan20}.(3)  presents a bijection between the sets  $\CI (A, \CD (A)-{\rm st.}, \CD (A)*A)$ and $ \CI (\CD (A), \CD (A)_{[0]})$.

\begin{proposition}\label{A3Jan20}
Let $A$ be an algebra. Then 
\begin{enumerate}
\item  $\Sub_{\CD (A)}(A)=\CI (A,  \CD (A)-{\rm st.})$.
\item  {\sc (Simplicity criterion for  ${}_{\CD (A)}A$)} The $\CD (A)$-module $A$ is simple iff there is no proper $\CD (A)$-stable ideal of $A$.
\item {\sc (The set $ \CI (\CD (A), \CD (A)_{[0]})$) } The map $$\CI (A, \CD (A)-{\rm st.}, \CD (A)*A)\ra \CI (\CD (A), \CD (A)_{[0]}), \;\;\ga \mapsto \ga + \CD (A)_{[0]}$$ is a bijection with the inverse $I\mapsto I\cap A$. The ideal $\CD (A)\CD (A)_{[0]}\CD (A)=\CD (A)*A+\CD (A)_{[0]}$ is the least ideal of the set $\CI (\CD (A), \CD (A)_{[0]})$, and the ideal $\CD (A)*A=A\cap \CD (A) \CD (A)_{[0]}\CD (A)$ is a $\CD (A)$-stable ideal of $A$.
\end{enumerate}
\end{proposition}

{\it Proof}. 1. Statement 1 is obvious.

2. Statement 2 follows from statement 1.

3. Let $\CD = \CD (A)$, $\CD_{[0]}=\CD (A)_{[0]}$ and $ \gb = \CD_{[0]}*A$.

(i) {\em If $I$ is an ideal of $\CD$ that contains $\CD_{[0]}$  then $I=\ga +\CD_{[0]}$ where $\ga =I\cap A$ is a $\CD$-stable ideal of $A$ such that $\gb\subseteq \ga$}: By (\ref{DA=ADA0}), $$I=I\cap \CD = I\cap (A+\CD_{[0]})=I\cap A+\CD_{[0]}=\ga +\CD_{[0]}.$$ 
 The left $\CD$-module $I/\CD_{[0]}=\ga$ is a submodule of the left $\CD$-module $A$. By statement 1, the ideal $\ga$ of $A$ is a $\CD$-stable ideal. Since $I\supseteq \CD\CD_{[0]}\CD = \gb +\CD_{[0]}$, $\ga \supseteq \gb$.
 
 (ii) {\em If $\ga$ is a $\CD$-stable ideal of $A$ that contains $\gb$ then $\ga+\CD_{[0]}$ is an ideal of $\CD$ that contains $\CD_{[0]}$}: $$\CD (\ga +\CD_{[0]})\CD = \CD \ga + \CD\CD_{[0]}\CD=\CD \ga + \gb + \CD_{[0]}=\CD*\ga+\gb+\CD_{[0]}\subseteq \ga +\CD_{[0]}.$$ Now, statement 3 follows from the statements (i) and (ii).  
$\Box $\\

Proposition \ref{B3Jan20} shows that if the algebra $A=\CA$ is a domain of essentially finite type over a perfect field then the ideal $\CD (\CA ) *\CA$ contains a power of the Jacobian ideal of the algebra $\CA$.

\begin{proposition}\label{B3Jan20}
Let $\CA$ be a domain of essentially finite type over a perfect field $K$ and $\ga_r$ be the Jacobian ideal of $\CA$. Then $\ga_r^i\subseteq \CD (\CA ) *\CA$ for some $i\geq 0$. 
\end{proposition}

{\it Proof}. Let $\CD = \CD (A)$, $\CD_{[0]}=\CD (A)_{[0]}$ and $ \gb = \CD_{[0]}*A$. By Proposition \ref{A2Jan20}.(3), the ideal $(\CD_{[0]})$ of $\CD$ is equal to $\gb+\CD_{[0]}$ and $\gb =\CA\cap (\CD_{[0]})$. By Theorem \ref{CXA9Jul19}.(3), $\ga_r^i\subseteq (\CD_{[0]})$ for some $i\geq 0$, and so $\gb \supseteq \CA \cap (\CD_{[0]})\supseteq \ga_r^i$.  $\Box $ \\

So, the subvariety $\Spec (\CA /\CD (\CA )*\CA)$ of $\Spec (\CA )$ consists of singular points of $\Spec (\CA )$.\\

If $\CA = K[x,y]/(y^2-x^3)$ is the algebra of regular functions on the cusp $y^2-x^3$ over a  field $K$ of characteristic zero then the algebra $\CD (\CA )$ is {\em simple},  \cite[Lemma 2.2.(2)]{SimCrit-difop}. Therefore, $\CD (\CA ) *\CA=\CA$ and $I=0$ in Theorem \ref{B3Jan20}. \\

{\bf The canonical form of a differential operator.}  For a finite set $\L$, we denote by $\N^{\L }$ the direct product  of $\L$ copies of the set of natural numbers $\N$. For an element  $\alpha = (\alpha_\l )$ of $\N^{\L }$, let $|\alpha |:=\sum_{\l \in \L} \alpha_\l$ and $(-1)^\alpha :=(-1)^{|\alpha |}$. For elements $\alpha , \beta \in \N^{\L }$, we write $\beta \leq \alpha$ if $\beta_\l \leq \alpha_\l$ for all $\l \in \L$. If $\beta \leq \alpha$ then ${\alpha\choose \beta}:=\prod_{\l \in \L}{\alpha_\l\choose \beta_\l}$.

\begin{theorem}\label{4Jan20}
Let $A$ be a finitely generated commutative algebra, $G=\{ x_\l\}_{\l \in \L}$ be a finite set of generators of $A$, $\CD (A)$ be the algebra of differential operators on $A$,   and $\{ \CD )A)_i\}_{i\geq 0}$ be the order filtration on $\CD (A)$.  
\begin{enumerate}
\item Each differential operator $\d \in \CD (A)_i$ of order $i$ is uniquely determined by the elements $\{ \ad^\alpha (\d )*1\, | \, \alpha \in \N^{\L }\}$ where $\ad^\alpha = \prod_{\l \in \L} \ad_{x_\l}^{\alpha_\l}$ for $\alpha = (\alpha_\l )\in \N^{\L }$.
\item For all elements $\alpha  \in \N^{\L }$ and $\d \in \CD (A)_i$, 
$$ \d (x^\alpha ) = \sum_{\beta \leq \alpha , |\beta |\leq i}(-1)^\beta {\alpha \choose \beta} \ad^\beta (\d )*1\cdot x^{\alpha-\beta}$$
where $x^{\alpha -\beta}=\prod_{\l \in \L} x_\l^{\alpha_\l - \beta_\l}$.
\end{enumerate}
\end{theorem}

{\it Proof}. 2. For the element $x_\l$,  we denote by $l_{x_\l}$ and $r_{x_\l}$ the left and right multiplication maps by the element $x_\l$, respectively. Now,
\begin{eqnarray*}
\d (x^\alpha) &=& \d x^\alpha*1=\d\prod_{\l \in \L}x_\l^{\alpha_\l}*1=\bigg( \prod_{\l \in \L}r_{x_\l}^{\alpha_\l}\d \bigg) *1=\bigg( \prod_{\l \in \L}(l_{x_\l}-\ad_{x_\l})^{\alpha_\l}\d \bigg) *1\\
 &=& \sum_{\beta \leq \alpha , |\beta |\leq i}(-1)^\beta {\alpha \choose \beta} \ad^\beta (\d )*1\cdot x^{\alpha-\beta}.
\end{eqnarray*}

1. Statement 1 follows from statement 2: Let  $\d , \d'\in \CD (A)_i$. If $\d = \d'$ then $\ad^\alpha (\d )*1=\ad^\alpha (\d' )*1$ for all elements $\alpha $ such that $|\alpha |\leq i$. The converse follows from statement 2.
 $\Box $\\



{\small

Department of Pure Mathematics

University of Sheffield

Hicks Building

Sheffield S3 7RH

UK

email: v.bavula@sheffield.ac.uk
}

\end{document}